# THE $G/GI/N$ QUEUE IN THE HALFIN–WHITT REGIME[1]

By Josh Reed

*New York University*

In this paper, we study the $G/GI/N$ queue in the Halfin–Whitt regime. Our first result is to obtain a deterministic fluid limit for the properly centered and scaled number of customers in the system which may be used to provide a first-order approximation to the queue length process. Our second result is to obtain a second-order stochastic approximation to the number of customers in the system in the Halfin–Whitt regime. This is accomplished by first centering the queue length process by its deterministic fluid limit and then normalizing by an appropriate factor. We then proceed to obtain an alternative but equivalent characterization of our limiting approximation which involves the renewal function associated with the service time distribution. This alternative characterization reduces to the diffusion process obtained by Halfin and Whitt [*Oper. Res.* **29** (1981) 567–588] in the case of exponentially distributed service times.

**1. Introduction.** In this paper, we study the $G/GI/N$ queue in the Halfin–Whitt regime. This problem has received considerable attention in the literature recently, however, to this date it has remained an open problem to extend the pioneering work of Halfin and Whitt [6] on the $GI/M/N$ queue to the more general $G/GI/N$ queue. In this paper and its sequel [18], we resolve this open problem by providing both fluid and diffusion limit results for the queue length process of the $G/GI/N$ queue in the Halfin–Whitt regime. In addition to providing these results, we also hope that the general methodology which is employed here, labeled the "Infinite Server Queue Systems Equations" approach (see below), will be helpful in future applications.

Received March 2007; revised April 2009.
[1]Supported in part by NSF Grant DMI-0300599.

Loosely speaking, the Halfin–Whitt regime is achieved by considering a sequence of many server queues indexed by the number of servers queues indexed by the number of servers $N$ where the arrival rate to the system grows large but the service time distribution remains fixed. Specifically, denoting by $\lambda^N$ the arrival rate to the $N$th system, we assume that

$$\lambda^N \to \infty \qquad \text{as } N \to \infty.$$

In order for the sequence of systems to remain stable, this then requires that the number of servers be large enough to handle the growing arrival rate. In particular, assuming without loss of generality that the mean service time is equal to one and denoting by $\rho^N = \lambda^N/N$ the traffic intensity to the $N$th system, in the Halfin–Whitt regime we assume that

$$\sqrt{N}(1 - \rho^N) \to \beta \qquad \text{as } N \to \infty,$$

where $-\infty < \beta < \infty$. Thus, in the Halfin–Whitt regime we assume that the traffic intensity of the system remains close to 1 while the number of servers grows without bound. Note also by the results of Kiefer and Wolfowitz [13], that, if $\beta > 0$, then, for large enough $N$, the sequence of systems will be stable in the pre-limit, while if $\beta < 0$, they will not. The case $\beta = 0$ is indeterminate.

Halfin and Whitt showed in their seminal paper [8] that in the regime described above, the properly centered and scaled queue length process will converge to a limiting diffusion. Unfortunately, they were not able to extend their results beyond the assumption of exponential service time distributions. This is mainly due to the fact that the infinitesimal generator approach to their proof breaks down when the service times are no longer exponentially distributed. This has naturally led to much speculation in the literature as to how to approach the situation of general service time distributions and specifically in such situations what the limiting process of the properly centered and scaled queue length process must be. In an effort to answer this question, several authors have recently obtained convergence results for carefully selected classes of service time distributions which are particularly well suited to analysis. Puhalskii and Reiman [17] have demonstrated convergence of the $G/PH/N$ queue length process in the Halfin and Whitt regime, where $PH$ stands for phase type service time distributions. Their approach is to consider a multi-dimensional Markovian process where each dimension corresponds to a different phase of the service time distribution. Jelenković, Mandelbaum and Momčilović [10] have shown convergence of the steady state distribution of the $GI/D/N$ queue, where $D$ stands for deterministic service times. Their proof involves focusing on a single server in the system and studying its queue length behavior as it evolves over time. Whitt in [21] has shown process level convergence of the $G/H_2^*/N/M$ queue, where $H_2^*$ stands for a mixture of an exponential random variable and a point mass at zero. In [16], Mandelbaum and Momčilović study the virtual waiting time process of $G/GI/N$ in the Halfin–Whitt regime assuming that the service time distribution possess finite support. Their approach relies on a combination of combinatorial and probabilistic arguments. Gamarnik and Momčilović [5] analyze the $GI/GI/N$ queue assuming that the service time distribution is lattice valued with finite support. They analyze the stationary values of the queue length and waiting time processes and show that the in the Halfin–Whitt regime, the diffusion scaled stationary value of the queue length converges to a limiting random variable corresponding to the stationary measure of a Markov chain, which, interestingly, may be recovered from our limit as well. Kaspi and Ramanan [12] consider service time distributions with a density and provide a fluid limit for the measure-valued process which keeps track of the amount of time that each customer has spent in the system. Nevertheless, with the exception of the results of Kapsi and Ramanan [12], it does not appear that any of the aforementioned approaches may be easily extended to the case of general service time distributions, and so to this date there has remained no general methodology for analyzing the $G/GI/N$ queue in the Halfin–Whitt regime. This is the main contribution of the present set of papers.

In particular, in this paper and its sequel [18], we provide two separate approaches for extending the results of Halfin and Whitt [6] to the $G/GI/N$ queue. Each of these approaches has its own unique set of advantages and disadvantages and in subsequent work we intend to provide important applications in which one approach may be more advantageous than the other. For the remainder of the present paper we focus our attention on the first approach which we label the "Infinite Server Queue System Equations" approach and defer discussion of the second approach, the "Idle Time System Equations" approach, until the sequel.

The main insight to the "Infinite Server Queue System Equations" approach is to write the system equations in a manner similar to the system equations for the $G/GI/\infty$ queue. Proposition 2.1 in Section 2 then provides a crucial link between our system and the $G/GI/\infty$ queue which allows the asymptotic analysis to proceed. In an effort to give a quick idea of what our main results, first recall the heavy traffic results found in Borovkov [2] and Krichagina and Puhalskii [14]. Recall that heavy traffic for the $G/GI/\infty$ queue is defined by letting the arrival rate to the system grow large while holding the service time distribution fixed. In such a regime, it can be shown that the properly centered and scaled queue length processes will converge to a Gaussian process which we denote by $\tilde{Q}_I$.

Let us therefore denote by $\tilde{Q}_I$ the limiting Gaussian process obtained for a $G/GI/\infty$ queue with the same sequence of arrival processes and an identical service time distribution as in our original sequence of $G/GI/N$ queues. Then the limiting process of Theorem 5.1 of Section 5 for the properly centered and scaled queue length process in our original sequence of $G/GI/N$ queues is given by the unique strong solution to

$$(1.1) \qquad \tilde{Q}(t) = \tilde{M}_Q(t) + \tilde{Q}_I(t) - \beta F_e(t) + \int_0^t \tilde{Q}^+(t-s)\, dF(s)$$

for $t \geq 0$, where $\tilde{Q}^+ = \max(\tilde{Q}, 0)$, $F$ is the CDF of the service time distribution, $F_e$ is the equilibrium distribution associated with $F$ [see (5.4)] and $\tilde{M}_Q$ is an additional process which is related to the initial conditions of the queue. Note that the additional integral term on the right-hand side of (1.1) is naturally positive as one would expect more customers in a $G/GI/N$ queue than in a corresponding $G/GI/\infty$ queue. Corollary 5.2 in Section 5 also shows that (1.1) may be equivalently expressed as

$$(1.2) \quad \tilde{Q}(t) = \zeta(t) + \int_0^t \zeta(t-s)\, dM(s) - \beta t - \int_0^t \tilde{Q}^-(t-s)\, dM(s)$$

for $t \geq 0$, where $\zeta = \tilde{M}_Q + \tilde{Q}_I$, $\tilde{Q}^- = \min(0, \tilde{Q})$ and $M$ is the renewal function associated with the pure renewal process with interarrival distribution $F$. From (1.2), it is then a matter of a few direct calculations to recover Halfin and Whitt's original results. We also point out that in Section 4 we develop fluid limit results which closely resemble (1.1) above.

The methodology of proof used in the present paper is heavily influenced by the results found in [14]. In particular, the authors in [14] use martingale techniques in order to show that certain processes associated with the queue length process are tight. In this paper, many of these same arguments are repeated again but with the slight modifications necessary in order to account for the finiteness of the number of servers. We therefore encourage the interested reader to concurrently review the results found in [14] in order to gain a fuller understanding of the present paper. In particular, one of the main insights from [14] is to show that the limiting process of $\tilde{Q}_I$ in (1.1) may be decomposed into a sum of two processes $\tilde{M}_1$ and $\tilde{M}_2$, which represent the randomness arising from the arrival process and service times, respectively. Furthermore, the process $\tilde{M}_2$ may represented as a double integral against the Kiefer process.

The remainder of this paper is now organized as follows. Section 2 provides the system equations for the $G/GI/N$ queue. In Section 3, we provide a regulator map result upon which our weak convergence argument will hinge. Sections 4 and 5 contain our weak convergence results for the queue length process. Specifically, in Section 4 we study the queue length process under fluid scaling and our main result in this section is Theorem 4.1. Next, in Section 5, we study the fluid centered queue length process under diffusion scaling in the Halfin–Whitt regime and our main result there is Theorem 5.1. Corollary 5.2 of Section 5 also provides an equivalent characterization of the limiting process obtained in Theorem 5.1. This then serves as the link between Halfin and Whitt's results and ours. In Section 6, we conclude by providing several directions for future research. The Appendix also includes several proofs which are similar in nature to those found in [14] but are necessary for our results and so are included here for completeness.

1.1. *Notation.* In what follows, all random variables are assumed to be defined on a common probability space $(\Omega, \mathcal{F}, \mathbb{P})$. Stochastic processes are assumed to measurable maps from $(\Omega, \mathcal{F})$ to $(D[0, \infty), \mathcal{D})$, where $D[0, \infty)$ is the space of all right continuous with left limit (RCLL) functions on $[0, \infty)$ and $\mathcal{D}$ is the Borel $\sigma$-algebra generated by the Skorohod $J_1$-topology, see Chapter 16 of [1] for further details.

We denote by $\mathcal{B}(\mathbb{R})$, the Borel $\sigma$-field on $\mathbb{R}$. For any two measure spaces, $(S_1, \mathcal{S}_1)$ and $(S_2, \mathcal{S}_2)$, we denote by $(S_1 \times S_2, \mathcal{S}_1 \times \mathcal{S}_2)$, the product measure space which is endowed with the product $\sigma$-field, $\mathcal{S}_1 \times \mathcal{S}_2$. Specifically, we define $(D^k[0, \infty), \mathcal{D}^k)$ to be the product measure space $(D[0, \infty) \times \cdots \times D[0, \infty), \mathcal{D} \times \cdots \times \mathcal{D})$.

We denote by $d_{J_1}$, the Skorohod metric on $D[0, \infty)$ and by $u$ the uniform metric. For each $x \in D[0, \infty)$ and $T \geq 0$, we denote by

$$\|x\|_T = \sup_{0 \leq t \leq T} |x(t)|,$$

the supremum metric on $[0, T]$. We also denote by $|\cdot|$, the Euclidian metric on $\mathbb{R}$. For any two metric spaces, $(S_1, m_1)$ and $(S_2, m_2)$, we denote by $(S_1 \times S_2, m_1 \times m_2)$, the product metric space which is endowed with the maximum metric $m_1 \times m_2$ defined by

$$(m_1 \times m_2)((x_1, x_2), (y_1, y_2)) = \max\{m_1(x_1, y_1), m_2(x_2, y_2)\}.$$

In particular, we define $(D^k[0, \infty), d^k_{J_1})$ to be the product metric space $(D^k[0, \infty), d^k_{J_1}) = (D[0, \infty) \times \cdots \times D[0, \infty), d_{J_1} \times \cdots \times d_{J_1})$ and we set $(D^k[0, \infty), u^k) = (D[0, \infty) \times \cdots \times D[0, \infty), u \times \cdots \times u)$.

**2. System equations for the $G/GI/N$ queue.** In this section, we provide the system equation for the $G/GI/N$ queue. One of the key insights from Halfin and Whitt [6] was that for large $N$, the $GI/M/N$ queue will, for stretches of time when the number of customers is low, behave as if it were an $GI/M/\infty$ queue. Our main results in Sections 4 and 5 show that the same holds true for the $G/GI/N$ queue as well. Our first step towards showing that this is the case is to write down the system equations for the $G/GI/N$ queue in a similar way to those for the $G/GI/\infty$ queue. For the reader's convenience, we will closely adhere to the notation used in [14] as many of the arguments we use here cite results from that paper.

Initially, at time $0-$, we assume that there are $Q_0$ customers in the system. The first $\min(Q_0, N)$ of these customers will be in service and the remainder are waiting to be served. Those customers in service at time $0-$ have already been in service for some amount of time and we denote by $\tilde{\eta}_i$ the residual service time of the $i$th customer in service at time $0-$. We assume that $\{\tilde{\eta}_i, i \geq 1\}$ is an i.i.d. sequence of random variables with common distribution $F_0$.

Customers next arrive to the system according to the arrival process $A = \{A(t), t \geq 0\}$ and are served on a first come first served (FCFS) basis. The arrival time of the $i$th customer is defined to be the quantity

$$\tau_i = \inf\{t \geq 0 : A(t) \geq i\}, \qquad i \geq 1.$$

Setting $\tau_0 = 0$, we also define

$$\xi_i = \tau_i - \tau_{i-1}, \qquad i \geq 1, \tag{2.1}$$

to be the interarrival times between the $(i-1)$st and $i$th customers to arrive to the system.

The $i$th customer to enter service after time $0-$ is assigned the service time $\eta_i$. We assume that $\{\eta_i, i \geq 1\}$ is an i.i.d. sequence of mean 1, random variables with common distribution $F$ whose tail distribution we denote by $G = 1 - F$. Note that we impose no assumptions on the service time distribution, other than it have a finite first moment.

For each $i \geq 1$, let $w_i$ denote the waiting time of the $i$th customer to arrive to the system after time $0-$ and let $\tilde{w}_i$ denote the waiting time of the $(N+i)$th initial customer in the system at time $0-$, if such a customer exists. We begin our indexing by $N+1$ since the first $N$ initial customers in the system will not have to wait. Using this notation as well as that of the previous paragraphs, the total number of customers in the system at time $t$ is given by

$$Q(t) = \sum_{i=1}^{\min(Q_0, N)} 1\{\tilde{\eta}_i > t\} + \sum_{i=1}^{(Q_0 - N)^+} 1\{\tilde{w}_i + \eta_i > t\} + \sum_{i=1}^{A(t)} 1\{\tau_i + w_i + \eta_{(Q_0 - N)^+ + i} > t\}. \tag{2.2}$$

We henceforth refer to the process $Q = \{Q(t), t \geq 0\}$ as the queue length process. It is important to note that $Q$ does not only count those customers in the queue waiting to be served but that indeed it counts the total number of customers in the system. The number of customers waiting to be served may however be recovered from $Q$ and is given by $(Q - N)^+$. Also note that in general $Q_0 \neq Q(0)$ since it is possible for customers to arrive to the system at time zero. One may think of $Q_0$ as being equal to $Q(0-)$, the left-hand limit of $Q$ at time $t = 0$.

Centering each of the indicator functions in the first summation on the right-hand side of (2.2) by the means and the indicator functions in the last two summations by their means conditional on their arrival times and waiting times, we obtain

(2.3)
$$Q(t) = \min(Q_0, N)\bar{F}_0(t) + W_0(t) + M_2(t)$$
$$+ \sum_{i=1}^{(Q_0-N)^+} G(t - \tilde{w}_i) + \sum_{i=1}^{A(t)} G(t - \tau_i - w_i),$$

where

(2.4)
$$W_0(t) = \sum_{i=1}^{\min(Q_0, N)} (1\{\tilde{\eta}_i > t\} - \bar{F}_0(t))$$

and

(2.5)
$$M_2(t) = \sum_{i=1}^{(Q_0-N)^+} (1\{\tilde{w}_i + \eta_i > t\} - G(t - \tilde{w}_i))$$
$$+ \sum_{i=1}^{A(t)} (1\{\tau_i + w_i + \eta_{(Q_0-N)^+ + i} > t\} - G(t - \tau_i - w_i)).$$

We also set $W_0 = \{W_0(t), t \geq 0\}$ and $M_2 = \{M_2(t), t \geq 0\}$.

Next, adding in and subtracting out the terms

(2.6)
$$A_G(t) = \int_0^t G(t - s) \, dA(s)$$

and $(Q_0 - N)^+ G(t)$ both to and from the right-hand side of (2.3), we obtain

(2.7)
$$Q(t) = I(t) + W_0(t) + M_2(t) + A_G(t)$$
$$+ \sum_{i=1}^{(Q_0-N)^+} (G(t - \tilde{w}_i) - G(t))$$
$$+ \sum_{i=1}^{A(t)} (G(t - \tau_i - w_i) - G(t - \tau_i)),$$

where

$$I(t) = \min(Q_0, N)\bar{F}_0(t) + (Q_0 - N)^+ G(t).$$

Note that $A_G(t)$ as given in (2.6) is the expected number of customers in $G/GI/\infty$ queue at time $t$ with same arrival process and service time distribution as in our $G/GI/N$ queue conditional on the arrival process $A$. We also set $A_G = \{A_G(t), t \geq 0\}$ and $I = \{I(t), t \geq 0\}$.

We now have the following key proposition.

PROPOSITION 2.1.  *For each $t \geq 0$,*

$$\sum_{i=1}^{A(t)} (G(t - \tau_i - w_i) - G(t - \tau_i))$$

$$= \int_0^t (Q(t-s) - N)^+ \, dF(s) - \sum_{i=1}^{(Q_0-N)^+} (G(t - \tilde{w}_i) - G(t)).$$

PROOF.  First note that for each time $t \geq 0$, we have that the total number of customers waiting to be served at time $t$ may be written as

$$(Q(t) - N)^+ = \sum_{i=1}^{(Q_0-N)^+} 1\{t < \tilde{w}_i\} + \sum_{i=1}^{A(t)} 1\{\tau_i \leq t < \tau_i + w_i\}.$$

We therefore have that

$$\sum_{i=1}^{A(t)} (G(t - \tau_i - w_i) - G(t - \tau_i))$$

$$= \sum_{i=1}^{A(t)} \int_{(t-(\tau_i+w_i))^+}^{t-\tau_i} dF(s)$$

$$= \sum_{i=1}^{A(t)} \int_0^\infty 1\{t - (\tau_i + w_i) < s \leq t - \tau_i\} \, dF(s)$$

$$= \sum_{i=1}^{A(t)} \int_0^\infty 1\{\tau_i \leq t - s < \tau_i + w_i\} \, dF(s)$$

$$= \int_0^\infty \sum_{i=1}^{A(t)} 1\{\tau_i \leq t - s < \tau_i + w_i\} \, dF(s)$$

$$= \int_0^t \left( (Q(t-s) - N)^+ - \sum_{i=1}^{(Q_0-N)^+} 1\{\tilde{w}_i > t - s\} \right) dF(s)$$

$$= \int_0^t (Q(t-s) - N)^+ \, dF(s) - \int_0^t \sum_{i=1}^{(Q_0-N)^+} 1\{\tilde{w}_i > t - s\} \, dF(s).$$

A reverse argument can now also be used to show that

$$\int_0^t \sum_{i=1}^{(Q_0-N)^+} 1\{\tilde{w}_i > t-s\}\, dF(s) = \sum_{i=1}^{(Q_0-N)^+} (G(t-\tilde{w}_i) - G(t)).$$

This completes the proof. □

Proposition 2.1 now allows us to rewrite equation (2.7) for the queue length at time $t$ as

(2.8)
$$Q(t) = I(t) + W_0(t) + M_2(t) + A_G(t)$$
$$+ \int_0^t (Q(t-s) - N)^+\, dF(s).$$

Equation (2.8) is the starting point for our analysis in Sections 4 and 5. In the next section, we develop a family of regulator map results which will be useful in representing the queue length process in (2.8).

**3. A family of regulator map results.** In this section, a family of regulator map results are provided which will be relied upon in the proof of our main results. In particular, these maps will provide convenient representations for the queue length processes in Sections 4 and 5.

Let $B$ be a cumulative distribution function on $\mathbb{R}$ and let $a \in \mathbb{R}$. For each $x \in D[0,\infty)$, we would like to find and characterize solutions $z \in D[0,\infty)$ to equations of the form

(3.1) $$z(t) = x(t) + \int_0^t (z(t-s) + a)^+\, dB(s), \qquad t \geq 0.$$

We therefore define the mapping $\varphi_B^a : D[0,\infty) \mapsto D[0,\infty)$ to be such that $\varphi_B^a(x)$ is a solution to (3.1) for each $x \in D[0,\infty)$. The following proposition now shows that $\varphi_B^a$ is uniquely defined and provides some regularity results for $\varphi_B^a$ as well. Its proof may be found in the Appendix.

PROPOSITION 3.1. *For each $x \in D[0,\infty)$, there exists a unique solution $\varphi_B^a(x)$ to (3.1). Moreover, the function $\varphi_B^a : D[0,\infty) \mapsto D[0,\infty)$ is Lipschitz continuous in the topology of uniform convergence over bounded intervals and measurable with respect to the Borel $\sigma$-field generated by the Skorohod $J_1$ topology.*

**4. Fluid limit results.** In this section, we obtain a nonlinear convolution equation as the fluid limit for the queue length process of the $G/GI/N$ queue in the Halin–Whitt regime. The limit which we obtain may be seen to be decomposed into four separate parts, one of which is the corresponding fluid limit for a $G/GI/\infty$ queue with the same sequence of arrival processes as our $G/GI/N$ queue and also with the same service time distribution. Although in many cases our fluid limit may not be directly solved for, we also present a special case in which it can which also highlights the rather unconventional behavior which our fluid limits may display.

Our underlying premise is that we are considering a sequence of $G/GI/N$ queues which we index by the number of servers $N$. In general, we will use a superscript $N$ to denote all processes and quantities associated with $N$th system.

Initially, at time $0-$, there are $Q_0^N$ customers in the $N$th system. The residual service time distribution of those customers in service in the $N$th system at time $0-$ are i.i.d. with common distribution $F_0$. We denote by $\{\tilde{\eta}_i, i \geq 1\}$ the i.i.d. sequence of residual service times.

Customers arrive to the $N$th system according to the arrival process $A^N = \{A^N(t), t \geq 0\}$. We denote by

$$\tau_i^N = \inf\{t \geq 0 : A^N(t) \geq i\}, \qquad i \geq 1,$$

the time of the arrival of the $i$th customer after time $0-$ to the $N$th system. The $i$th customer to enter service after time $0-$ in the $N$th system is assigned the service $\eta_i$, where $\{\eta_i, i \geq 1\}$ is an i.i.d. sequence of random variables with common distribution $F$. Finally, we denote by $G = 1 - F$ the tail distribution of $F$. Note that neither the sequence of residual service times or actual service times is changing with $N$.

For the remainder of this section, we study the fluid scaled queue length process $\bar{Q}^N = \{N^{-1}Q^N(t), t \geq 0\}$. Our high level approach will be similar in spirit to fluid limit proofs for the queue length process in conventional heavy-traffic in which the number of servers remains fixed but the service rate is increased. In particular, we first provide a representation of the queue length process in terms of the regulator mapping $\varphi_F^a$ provided by Proposition 3.1 and an associated free process, say $\bar{X}$. We then provide several weak convergence results related to $\bar{X}$ which may be used in conjunction with the Continuous Mapping theorem and the representation in terms of $\varphi_F^a$ in order to establish the main result of the section, Theorem 4.1, which details the asymptotic behavior of the fluid scaled queue length process.

Let $Q^N = \{Q^N(t), t \geq 0\}$ be the queue length process in the $N$th system and recall that by equation (2.8) of Section 2, we have that

$$Q^N(t) = I^N(t) + W_0^N(t) + M_2^N(t) + A_G^N(t)$$
(4.1)
$$+ \int_0^t (Q^N(t-s) - N)^+ \, dF(s).$$

If we now define the fluid scaled quantities,

(4.2) $$\bar{Q}^N(t) = \frac{Q^N(t)}{N},$$

$$\bar{I}^N(t) = \frac{I^N(t)}{N},$$

$$\bar{W}_0^N(t) = \frac{W_0^N(t)}{N},$$

(4.3) $$\bar{M}_2^N(t) = \frac{M_2^N(t)}{N}$$

and

(4.4) $$\bar{A}_G^N(t) = \frac{A_G^N(t)}{N},$$

it then follows from (4.1) that

(4.5) $$\bar{Q}^N(t) = \bar{I}^N(t) + \bar{W}_0^N(t) + \bar{M}_2^N(t) + \bar{A}_G^N(t)$$
$$+ \int_0^t (\bar{Q}^N(t-s) - 1)^+ \, dF(s).$$

Furthermore, since by Proposition 3.1, the mapping $\varphi_F^a$ is uniquely defined with $a = -1$, setting

$$\bar{Q}^N = \{\bar{Q}^N(t), t \geq 0\},$$
$$\bar{I}^N = \{\bar{I}^N(t), t \geq 0\},$$
$$\bar{W}_0^N = \{\bar{W}_0^N(t), t \geq 0\},$$
$$\bar{M}_2^N = \{\bar{M}_2^N(t), t \geq 0\}$$

and

$$\bar{A}_G^N(t) = \{\bar{A}_G^N(t), t \geq 0\},$$

we have from (4.5) that

(4.6) $$\bar{Q}^N = \varphi_F^a(\bar{I}^N + \bar{W}_0^N + \bar{M}_2^N + \bar{A}_G^N),$$

with $a = -1$. The representation (4.6) above will turn out to be useful when proving the main result of this section.

We now state several preliminary results in preparation for the statement of the main result of the section, Theorem 4.1. Our first result shows that $\bar{W}_0^N$ converges to zero as $N$ goes to $\infty$.

PROPOSITION 4.1. $\bar{W}_0^N \Rightarrow 0$ as $N \to \infty$.

PROOF. First, note that

$$\bar{W}_0^N(t) = N^{-1} \sum_{i=1}^{N \min(N^{-1} Q_0^N, 1)} (1\{\tilde{\eta}_i > t\} - \bar{F}_0(t)).$$

Thus, for each $T > 0$ and $\delta > 0$, we have

$$P\left(\sup_{0 \leq t \leq T}\left|N^{-1}\sum_{i=1}^{N\min(N^{-1}Q_0^N, 1)}(1\{\tilde{\eta}_i > t\} - \bar{F}_0(t))\right| > \delta\right)$$

$$\leq P\left(\sup_{0 \leq x \leq 1}\sup_{0 \leq t \leq T}\left|N^{-1}\sum_{i=1}^{\lfloor xN \rfloor}(1\{\tilde{\eta}_i > t\} - \bar{F}_0(t))\right| > \delta\right).$$

However, by Lemma 3.1 in [14],

$$P\left(\sup_{0 \leq x \leq 1}\sup_{0 \leq t \leq T}\left|N^{-1}\sum_{i=1}^{\lfloor xN \rfloor}(1\{\tilde{\eta}_i > t\} - \bar{F}_0(t))\right| > \delta\right) \to 0 \quad \text{as } N \to \infty,$$

which completes the proof. $\square$

We next show that $\bar{M}_2^N$ converges in distribution to zero. The full proof of this result may be found in the Appendix.

PROPOSITION 4.2. $\bar{M}_2^N \Rightarrow 0$ as $N \to \infty$.

PROOF. See the Appendix. $\square$

The following is now the main result of this section. It provides a deterministic first-order approximation to the queue length process. Later, in Section 5, we use this result to center the queue length process and obtain a second-order approximation.

Let

$$\bar{Q}_0^N = \frac{Q_0^N}{N}$$

be the fluid scaled initial number of customer in the system at time $0-$ and

$$\bar{A}^N(t) = \frac{A^N(t)}{N}$$

be the fluid scaled number of arrivals by time $t \geq 0$. We also set $\bar{A}^N = \{\bar{A}^N(t), t \geq 0\}$ to be the fluid scaled arrival process. We then have the following.

THEOREM 4.1. If $(\bar{Q}_0^N, \bar{A}^N) \Rightarrow (\bar{Q}_0, \bar{A})$ in $(\mathbb{R} \times D[0, \infty), |\cdot| \times d_{J_1})$ as $N \to \infty$, where $\bar{A}$ is a stochastic process with $\mathbb{P}$-a.s. continuous sample paths, then $\bar{Q}^N \Rightarrow \bar{Q}$ as $N \to \infty$, where $\bar{Q}$ is the unique strong solution to

$$\bar{Q}(t) = \min(\bar{Q}_0, 1)\bar{F}_0(t) + (\bar{Q}_0 - 1)^+ G(t)$$
(4.7)
$$+ \int_0^t G(t-s)\,d\bar{A}(s) + \int_0^t (\bar{Q}(t-s) - 1)^+\,dF(s)$$

for $t \geq 0$.

PROOF. First, note that by the definition of $\bar{A}_G^N$ in (4.4) and the assumption of the theorem that $\bar{A}^N \Rightarrow \bar{A}$ as $N \to \infty$, where $\bar{A}$ has $\mathbb{P}$-a.s. continuous sample paths, it follows as in the proof of Theorem 3 of [14] that

$$
\begin{aligned}
\bar{A}_G^N &= \int_0^{\cdot} G(\cdot - s) \, d\bar{A}^N(s) \\
&\Rightarrow \int_0^{\cdot} G(\cdot - s) \, d\bar{A}(s)
\end{aligned}
\tag{4.8}
$$

as $N \to \infty$. Next, setting

$$\bar{M}_3^N = \bar{W}_0^N + \bar{M}_2^N + \bar{A}_G^N,$$

it follows by Propositions 4.1 and 4.2 and (4.8) that

$$\bar{M}_3^N \Rightarrow \int_0^{\cdot} G(\cdot - s) \, d\bar{A}(s) \qquad \text{as } N \to \infty.$$

Since, by assumption, $(\bar{Q}_0^N, \bar{A}^N) \Rightarrow (\bar{Q}_0, \bar{A})$ in $(\mathbb{R} \times D[0, \infty), |\cdot| \times d_{J_1})$ as $N \to \infty$, it now follows by Theorem 11.4.5 in [20] that

$$(\bar{M}_3^N, \bar{Q}_0^N) \Rightarrow \left( \int_0^{\cdot} G(\cdot - s) \, d\bar{A}(s), \bar{Q}_0 \right) \qquad \text{in } (\mathbb{R} \times D[0, \infty), |\cdot| \times d_{J_1})$$

as $N \to \infty$. By Theorem 11.4.1 in [20], the space $\mathbb{R} \times D$ is separable under the product topology induced by the maximum metric $|\cdot| \times d_{J_1}$ and thus, by the Skorohod representation theorem [20], there exists some alternate probability space, $(\hat{\Omega}, \hat{\mathcal{F}}, \hat{P})$, on which are defined a sequence of processes

$$\{(\hat{M}_3^N, \hat{Q}_0^N), N \geq 1\} \tag{4.9}$$

such that

$$(\hat{M}_3^N, \hat{Q}_0^N) \stackrel{d}{=} (\bar{M}_3^N, \bar{Q}_0^N) \qquad \text{for } N \geq 1, \tag{4.10}$$

and also processes

$$\left( \int_0^{\cdot} G(\cdot - s) \, d\bar{A}(s), \hat{Q}_0 \right) \stackrel{d}{=} \left( \int_0^{\cdot} G(\cdot - s) \, d\bar{A}(s), \bar{Q}_0 \right), \tag{4.11}$$

where

$$(\hat{M}_3^N, \hat{Q}_0^N) \to \left( \int_0^{\cdot} G(\cdot - s) \, d\bar{A}(s), \hat{Q}_0 \right)$$

$$\text{in } (\mathbb{R} \times D, |\cdot| \times d_{J_1}) \; \hat{\mathbb{P}}\text{-a.s.} \tag{4.12}$$

as $N \to \infty$. Furthermore, as the process $\int_0^{\cdot} G(\cdot - s) \, d\bar{A}(s)$ on the right-hand side of (4.12) is, by the assumption of the continuity on $\bar{A}$, continuous, it follows that the convergence in (4.12) can also be strengthened to convergence in $(\mathbb{R} \times D, |\cdot| \times u)$.

Now set
$$\hat{I}^N = \min(\hat{Q}_0^N, 1)\bar{F}_0 + (\hat{Q}_0^N - 1)^+ G$$
and note that by (4.10), we have

(4.13) $\quad (\hat{M}_3^N, \hat{I}^N) \stackrel{d}{=} (\bar{M}_3^N, \bar{I}^N) \quad$ for $N \geq 1$.

Furthermore, letting
$$\hat{I} = \min(\hat{Q}_0, 1)\bar{F}_0 + (\hat{Q}_0 - 1)^+ G,$$
we have for each $T \geq 0$, by (4.12),

$$
\begin{aligned}
(4.14) \quad \sup_{0 \leq t \leq T} |\hat{I}^N(t) - \hat{I}(t)| &= \sup_{0 \leq t \leq T} |(\min(\hat{Q}_0^N, 1) - \min(\hat{Q}_0, 1))\bar{F}_0(t) \\
&\quad + ((\hat{Q}_0^N - 1)^+ - (\hat{Q}_0 - 1)^+)G(t)| \\
&\leq |\hat{Q}_0^N - \hat{Q}_0| \sup_{0 \leq t \leq T} (\bar{F}_0(t) + G(t)) \\
&\leq 2|\hat{Q}_0^N - \hat{Q}_0|
\end{aligned}
$$

(4.15) $\quad\quad\quad\quad\quad\quad\quad \to 0 \quad \hat{\mathbb{P}}\text{-a.s. as } N \to \infty.$

Now let
$$\hat{Q}^N = \varphi_F^a(\hat{I}^N + \hat{M}_3^N),$$
where $a = -1$, and note that by the representation (4.6), (4.9), (4.13) and the measurability of $\varphi_F^a$ from Proposition 3.1, it follows that

(4.16) $\quad\quad\quad\quad \hat{Q}^N \stackrel{d}{=} \tilde{Q}^N \quad$ for $N \geq 1$.

Furthermore, it follows from (4.12), (4.14) and the continuity of $\varphi_F^a$ with respect to the topology of uniform convergence over compact sets from Proposition 3.1, that

$$
\begin{aligned}
\hat{Q}^N &= \varphi_F^a(\hat{I}^N + \hat{M}_3^N) \\
&\to \varphi_F^a\left(\min(\hat{Q}_0, 1)\bar{F}_0 + (\hat{Q}_0 - 1)^+ G + \int_0^\cdot G(\cdot - s)\,d\bar{A}(s)\right)
\end{aligned}
$$

in $(D[0, \infty), u)$ $\hat{\mathbb{P}}$-a.s. as $N \to \infty$. Thus, since convergence in $(D, u)$ implies convergence in $(D, d_{J_1})$ and almost sure convergence implies convergence in distribution, it follows by the measurability of $\psi_F^a : (D[0, \infty), \mathcal{D}) \mapsto (D[0, \infty), \mathcal{D})$ from Proposition 3.1 and (4.16) that

$$\tilde{Q}^N \Rightarrow \varphi_F^a\left(\min(\bar{Q}_0, 1)\bar{F}_0 + (\bar{Q}_0 - 1)^+ G + \int_0^\cdot G(\cdot - s)\,d\bar{A}(s)\right) \quad \text{as } N \to \infty,$$

which completes the proof. $\square$

Note that the fluid limit $\bar{Q}$ given by (4.7) of Theorem 4.1 may be decomposed into four separate parts. The first two terms on the right-hand side of (4.7) are representative of the fluid scaled number of customers in the queue at time $0-$. Specifically, $\min(\bar{Q}_0, 1)\bar{F}_0(t)$ is the limiting fluid scaled number of customers who were in the service at time $0-$ and still remain in the system at time $t$ and $(\bar{Q}_0 - 1)^+ G(t)$ is representative of those customers who were waiting in the queue to be served at time $0-$. Next, the term $\int_0^t G(t-s) \, d\bar{A}(s)$ may be viewed as the limiting fluid scaled number of customers in the system at time $t$ in an $G/GI/\infty$ queue with the same sequence of arrival processes and service time distribution as in our $G/GI/N$ queue and which starts out empty at time $0-$. Finally, the integral term on the right-hand side of (4.7) may be thought of as an adjustment to the infinite server term immediately preceding it which takes into account the waiting times of customers.

In general, the limiting process of Theorem 5.1, $\bar{Q}$, cannot be directly solved for. This is mainly due to the presence of the nonlinear $()^+$ operator in the integral term. However, under certain special circumstances it can. The following example now presents one such case in which an explicit solution may be found.

EXAMPLE 1. Consider the case of deterministic service times in which we have $\bar{Q}_0 = 1$ and $\bar{A} = e$. Further, we also assume that the residual service times are constant with mean equal to 1 so that $F_0(x) = F(x) = 1\{x \geq 1\}, x \geq 0$. In this case, (4.7) takes the rather simple form

$$\bar{Q}(t) = 1 + t, \qquad 0 \leq t < 1, \tag{4.17}$$

and

$$\bar{Q}(t) = 1 + (\bar{Q}(t) - 1)^+, \qquad t \geq 1. \tag{4.18}$$

Solving this recursion, one finds that

$$\bar{Q}(t) = 1 + t - \lfloor t \rfloor, \qquad t \geq 0. \tag{4.19}$$

Thus, $\bar{Q}$ exhibits the sawtooth pattern as shown in Figure 1 above. Note the rather unconventional nature of the fluid limit in Figure 1. In particular, it is periodic with a period of 1 and it is also discontinuous. Thus, our limit process may in general display rather irregular behavior. However, as is shown in Section 5 below, if one starts out the queue length process under general "equilibrium" conditions, then much more regular behavior of the fluid limit may be obtained.

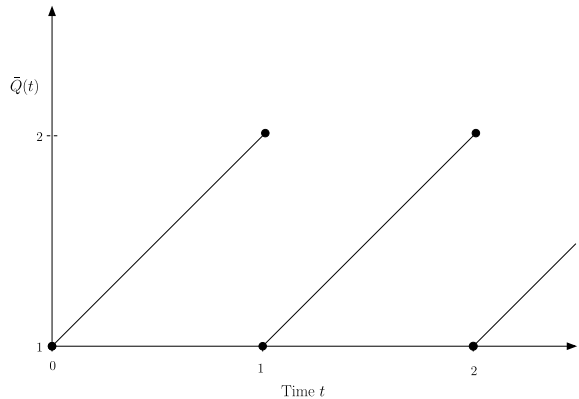

Fig. 1. *The graph of $\{\bar{Q}(t), t \geq 0\}$ for $\bar{Q}_0 = 1, \bar{A} = e$ and $F_0(t) = F(t) = 1\{t \geq 1\}$.*

**5. Diffusion limit results.** In this section, we obtain limiting results for the $G/GI/N$ queue in the Halfin–Whitt regime. Our main result of the section is to provide a limiting approximation for the diffusion scaled queue length process in this regime. This limiting approximation may be viewed as the solution to a stochastic nonlinear convolution equation. In order to proceed, we first express the queue length process for the $G/GI/N$ queue via the regulator map $\varphi_F^a$ defined in Section 3. We then provide several useful propositions in preparation for the statement of our main result, Theorem 5.1, which provides a limiting approximation for the diffusion scaled queue length process. In Corollary 5.2, we provide an alternative representation of the limiting process obtained in Theorem 5.1. This representation has several desirable properties and we conclude the section by showing how in the case of exponentially distributed service times and renewal arrivals it reduces to the diffusion obtained by Halfin and Whitt [7].

5.1. *The Halfin–Whitt regime heavy-traffic regime.* In order to proceed, we must first provide a detailed description of the Halfin–Whitt regime. Our setup is similar to Section 4 where we obtained our fluid limit results. Again, we consider a sequence of $G/GI/N$ queues indexed by the number of servers $N$. Initially, at time $0-$, there are $Q_0^N$ customers in the system and the first $\min(Q_0^N, N)$ of these customers have i.i.d. residual service times with common distribution $F_0$. We denote by $\tilde{\eta}_i$, the residual service time of the $i$th customer in service at time $0-$.

Customers arrive to the $N$th system according to the arrival process $A^N = \{A^N(t), t \geq 0\}$ and we denote by

$$\tau_i^N = \inf\{t \geq 0 : A^N(t) \geq i\}, \qquad i \geq 1,$$

the time of the arrival of the $i$th customer after time $0-$ to the $N$th system. We also assume that there exists a sequence of constants $\{\rho^N, N \geq 1\}$ such that $\rho^n \to 1$ as $N \to \infty$, and, where, setting

(5.1) $$\tilde{A}^N(t) = \frac{A^N(t) - N\rho^N t}{\sqrt{N}}, \qquad t \geq 0,$$

and $\tilde{A}^N = \{\tilde{A}^N(t), t \geq 0\}$, we have that

(5.2) $$\tilde{A}^N \Rightarrow \tilde{\xi} \qquad \text{as } N \to \infty,$$

where $\tilde{\xi}$ is a stochastic process with $\mathbb{P}$-a.s. continuous sample paths. Loosely speaking, one may interpret $\rho^N$ as the arrival rate of customers to the $N$th system.

Note that assumption (5.2) is flexible from a modeling point of view. In heavy traffic theory, one often assumes that $A^N(e) = A(Ne)$, where $A$ is a renewal process, in which case, by Donsker's theorem, the process $\tilde{\xi}$ in (5.2) turns out to be a Brownian motion. The interpretation this assumption is that customers are emanating from a single source, albeit at a rapid rate. However, in many applications, with telephone call centers being just one such example, it is perhaps more natural to assume that customers are emanating from many sources. This then leads to the assumption that $A^N$ is a superposition of many i.i.d. renewal arrival processes, that is, $A^N = \sum_{i=1}^{N} A_i$. Under such an assumption, the process $\tilde{\xi}$ turns out to be a centered Gaussian process whose covariance structure is inherited from that of each of the individual $A_i$'s. The interested reader is referred to Section 7.2 of [20] for further details on this remark.

As in Section 4, the service time distribution is held fixed as we index through $N$. We therefore denote by $\eta_i$ the service time of the $i$th customer to enter service after time $0-$, where $\{\eta_i, i \geq 1\}$ is an i.i.d. sequence of mean 1 random variables with common distribution $F$. We denote by $G = 1 - F$ the tail distribution of $F$. Note also that we place no additional restrictions on $F$ beyond a first moment.

It now remains to provide the key relationship characterizing the Halfin–Whitt regime. As noted above, we have that by (5.1) and (5.2), the quantity $N\rho^N$ may be loosely interpreted as the arrival rate to the $N$th system. Next, since there are $N$ servers in the $N$th system and it is assumed that the service rate is fixed at one, it also follows that $\rho^N$ may also be interpreted as the traffic intensity of the $N$th system. The Halfin–Whitt regime is now achieved by specifying the rate at which the traffic intensity of the system converges to one as $N$ grows to infinity. Specifically, we assume that

(5.3) $$\sqrt{N}(1 - \rho^N) \to \beta \qquad \text{as } N \to \infty,$$

where $-\infty < \beta < \infty$.

5.2. *Initial conditions.* In proving our main diffusion limit result, it will be useful to assume that the limiting fluid scaled number of customers in the system is constant for all $t \geq 0$. Let

$$(5.4) \qquad F_e(x) = \int_0^x G(u)\,du, \qquad x \geq 0,$$

be the equilibrium distribution associated with $F$. The following result may now be seen as a corollary to Theorem 5.1 of Section 4.

COROLLARY 5.1.  *If $F_0 = F_e$ and $\bar{Q}_0^N \Rightarrow 1$ as $N \to \infty$, then $\bar{Q}^N \Rightarrow 1$ as $N \to \infty$.*

PROOF. First note that (5.1), (5.2) and the Halfin–Whitt assumption (5.3) imply the functional weak large of law large numbers result, $\bar{A}^N \Rightarrow e$ as $N \to \infty$, where $e = \{t, t \geq 0\}$ is the identity process. Thus, by Theorem 3.9 in [1] and the assumption $\bar{Q}_0^N \Rightarrow 1$ as $N \to \infty$, we have that $(\bar{Q}_0^N, \bar{A}^N) \Rightarrow (1, e)$ in $(\mathbb{R} \times D, |\cdot| \times d_{J_1})$ as $N \to \infty$.

It now follows by Theorem 5.1 in Section 4 that $\bar{Q}^N \Rightarrow \bar{Q}$ as $N \to \infty$, where $\bar{Q}$ is given by the unique solution to

$$\begin{aligned}\bar{Q}(t) &= \min(1,1)\bar{F}_e(t) + (1-1)^+ G(t) + \int_0^t G(t-s)\,de(s) \\ &\quad + \int_0^t (\bar{Q}(t-s) - 1)^+\,dF(s) \\ &= \bar{F}_e(t) + \int_0^t G(t-s)\,ds + \int_0^t (\bar{Q}(t-s) - 1)^+\,dF(s) \\ &= 1 + \int_0^t (\bar{Q}(t-s) - 1)^+\,dF(s)\end{aligned}$$

for $t \geq 0$. By inspection, one sees that $\bar{Q}(t) = 1$ for $t \geq 0$ is the unique solution to this equation, which completes the proof.  □

For the remainder of this section, we assume that the simplifying assumptions of the above corollary hold. That is, we assume that the initial residual service time distribution is equal to the equilibrium distribution $F_e$ and that $\bar{Q}^N \Rightarrow 1$ as $N \to \infty$. In future papers, we intend to remove these assumptions.

5.3. *Weak convergence results.* We now proceed to provide our weak convergence results. We begin by providing a convenient representation for the queue length process in terms of the regulator map of Proposition 3.1.

Recall that by equation (2.8) of Section 2, we have that the queue length at time $t$ may be written as

(5.5)
$$Q^N(t) = I^N(t) + W_0^N(t) + M_2^N(t) + A_G^N(t)$$
$$+ \int_0^t (Q^N(t-s) - N)^+ \, dF(s).$$

We would next like to the center the queue length process by its fluid limit, appropriately scaled. By Corollary 5.1 above, we have that the limiting fluid number of customers in the system, $\bar{Q}$, is equal to 1 for all $t \geq 0$. We therefore choose to center the queue length process in the $N$th system by $N$. Performing such a centering as well as some algebraic manipulations and recalling the definition of $F_e$ from (5.4), one then obtains that

(5.6)
$$Q^N(t) - N = M_Q^N(t) + H^N(t) + W_0^N(t) + M_2^N(t) + M_1^N(t)$$
$$- N(1 - \rho^N) F_e(t) + \int_0^t (Q^N(t-s) - N)^+ \, dF(s),$$

where

$$M_Q^N(t) = (Q_0^N - N)^+ (G(t) - \bar{F}_e(t)),$$
$$H^N(t) = (Q_0^N - N) \bar{F}_e(t)$$

and

(5.7)
$$M_1^N(t) = \int_0^t G(t-s) \, d(A^N(s) - N\rho^N s).$$

Let $M_Q^N = \{M_Q^N(t), t \geq 0\}, H^N = \{H^N(t), t \geq 0\}$ and $M_1^N = \{M_1^N(t), t \geq 0\}$.

If we now define the diffusion scaled quantities,

(5.8)
$$\tilde{Q}^N(t) = \frac{Q^N(t) - N}{\sqrt{N}},$$
$$\tilde{M}_Q^N(t) = \frac{M_Q^N(t)}{\sqrt{N}},$$
$$\tilde{H}^N(t) = \frac{H^N(t)}{\sqrt{N}},$$
$$\tilde{W}_0^N(t) = \frac{W_0^N(t)}{\sqrt{N}},$$
$$\tilde{M}_1^N(t) = \frac{M_2^N(t)}{\sqrt{N}}$$

and

(5.9)
$$\tilde{M}_2^N(t) = \frac{M_1^N(t)}{\sqrt{N}},$$

it then follows, dividing (5.6) by $\sqrt{N}$, that

$$\tilde{Q}^N(t) = \tilde{M}_Q^N(t) + \tilde{Q}_I^N(t) - \sqrt{N}(1-\rho^N)F_e(t) \tag{5.10}$$
$$+ \int_0^t \tilde{Q}^{N,+}(t-s)\,dF(s),$$

where

$$\tilde{Q}_I^N(t) = \tilde{H}^N(t) + \tilde{W}_0^N(t) + \tilde{M}_1^N(t) + \tilde{M}_2^N(t).$$

Letting

$$\tilde{Q}^N = \{\tilde{Q}^N(t), t \geq 0\},$$
$$\tilde{M}_Q^N = \{\tilde{M}_Q^N(t), t \geq 0\},$$
$$\tilde{H}^N = \{\tilde{H}^N(t), t \geq 0\},$$
$$\tilde{W}_0^N = \{\tilde{W}_0^N(t), t \geq 0\},$$
$$\tilde{M}_2^N = \{\tilde{M}_2^N(t), t \geq 0\},$$
$$\tilde{M}_1^N = \{\tilde{M}_1^N(t), t \geq 0\}$$

and

$$\tilde{Q}_I^N = \tilde{H}^N + \tilde{W}_0^N + \tilde{M}_1^N + \tilde{M}_2^N,$$

we then have, since the mapping $\varphi_F^a$ with $a=0$ is by Proposition 3.1 uniquely defined, that (5.10) may also be written as

$$\tilde{Q}^N = \varphi_F^a(\tilde{M}_Q^N + \tilde{Q}_I^N - \sqrt{N}(1-\rho^N)F_e), \tag{5.11}$$

with $a=0$.

The representation (5.11) will be useful in the proof our main result. However, before stating this result, we first provide several preliminary propositions and lemmas which are interesting in their own right and will be crucial in the proof of our main result.

The proof of the following result may now be found in the Appendix. It provides a limiting process which represents the randomness of the service times in the limit. For further details on the limit process below, the interested reader may consult [14].

PROPOSITION 5.1. *Let*

$$\hat{M}_2^N(t) = N^{-1/2} \sum_{i=1}^{Nt} (1\{N^{-1}i + \eta_i \geq t\} - G(t - N^{-1}i)), \qquad t \geq 0,$$

*and set* $\hat{M}_2^N = \{\hat{M}_2^N(t), t \geq 0\}$. *Then,*

$$(\tilde{M}_2^N, \hat{M}_2^N) \Rightarrow (\tilde{M}_2, \tilde{M}_2) \qquad \text{in } (D^2[0,\infty), d_{J_1}^2) \text{ as } N \to \infty,$$

where $\tilde{M}_2$ is a centered Gaussian process with covariance structure

$$E[\tilde{M}_2(t)\tilde{M}_2(t+\delta)] = \int_0^t G(t+\delta-u)F(t-u)\,du \qquad \text{for } t, \delta \geq 0.$$

PROOF. See the Appendix. □

We next prove a joint convergence result on the diffusion scaled processes defined at the beginning of this section. Let

$$\tilde{Q}_0^N = \frac{Q_0^N - N}{\sqrt{N}} \tag{5.12}$$

be the diffusion scaled number of customers in the system at time $0-$.

Next, let $\tilde{W}_0 = \{\tilde{W}_0(t), t \geq 0\}$ be a Brownian bridge. In other words, $\tilde{W}_0$ is the unique continuous, centered Gaussian process on $[0,1]$ with covariance function

$$E[\tilde{W}_0(s)\tilde{W}_0(t)] = (s \wedge t) - st, \qquad 0 \leq s \leq t \leq 1.$$

Moreover, set $\hat{W}_0(F_e) = \{\tilde{W}_0(F_e(t)), t \geq 0\}$, where $F_e$ is the equilibrium distribution associated with $F$ as defined in (5.4). One may view $\tilde{W}_0(F_e)$ as a time changed Brownian bridge.

We then have the following result.

PROPOSITION 5.2. *If $\tilde{Q}_0^N \Rightarrow \tilde{Q}_0$ as $N \to \infty$, then*

$$(\tilde{Q}_0^N, \tilde{W}_0^N, \tilde{A}^N, \tilde{M}_2^N) \Rightarrow (\tilde{Q}_0, \tilde{W}_0(F_e), \tilde{\xi}, \tilde{M}_2) \qquad \text{in } (\mathbb{R} \times D^3[0,\infty), |\cdot| \times d_{J_1}^3)$$

*as $N \to \infty$, where each of the limiting processes appearing on the right-hand side above are independent of one another.*

PROOF. We first show convergence of the marginals. The convergence of $\tilde{Q}_0^N$ to $\tilde{Q}_0$ is clear by assumption (5.12).

Let

$$\hat{W}^N(t) = N^{-1/2} \sum_{i=1}^N (1\{\tilde{\eta}_i > t\} - \bar{F}_e(t))$$

and set $\hat{W}^N = \{\hat{W}^N(t), t \geq 0\}$. The convergence

$$(5.13) \quad (\hat{W}_0^N, \tilde{W}_0^N) \Rightarrow (\tilde{W}_0(F_e), \tilde{W}_0(F_e)) \qquad \text{in } (D^2[0,\infty), d_{J_1}^2) \text{ as } N \to \infty,$$

follows by the representation

$$\tilde{W}^N(t) = N^{-1/2} \sum_{i=1}^{N\min(N^{-1}Q_0^N, 1)} (1\{\tilde{\eta}_i > t\} - \bar{F}_e(t)),$$

the Random Time Change theorem [1] and Lemma 3.1 of [14], since, by the Continuous Mapping theorem and assumption (5.12),

$$\min(N^{-1}Q_0^N, 1) \Rightarrow 1 \qquad \text{as } N \to \infty.$$

The convergence of $\tilde{A}^N$ to $\tilde{\xi}$ follows by assumption (5.2) and the convergence of $\tilde{M}_2^N$ to $\tilde{M}_2$ is immediate by Proposition 5.1.

It remains to show the joint convergence as stated in the proposition. The convergence

$$(\tilde{Q}_0^N, \hat{W}_0^N, \tilde{A}^N, \hat{M}_2^N) \Rightarrow (\tilde{Q}_0, \tilde{W}_0(F_e), \tilde{\xi}, \tilde{M}_2) \qquad \text{in } (\mathbb{R} \times D^3[0, \infty), |\cdot| \times d_{J_1}^3)$$

as $N \to \infty$, follows by Theorem 11.4.4 in [20] since each of the component processes appearing in the prelimit above are independent of one another and further, they converge to their desired limits as shown in the previous paragraph. Next, note that

$$(|\cdot| \times d_{J_1}^3)((\tilde{Q}_0^N, \tilde{W}_0^N, \check{A}^N, \tilde{M}_2^N), (\tilde{Q}_0^N, \hat{W}_0^N, \check{A}^N, \hat{M}_2^N))$$
$$\leq d_{J_1}(\tilde{W}_0^N, \hat{W}_0^N) + d_{J_1}(\tilde{M}_2^N, \hat{M}_2^N)$$

and thus, if we can show that

(5.14) $\qquad d_{J_1}(\tilde{W}_0^N, \hat{W}_0^N) + d_{J_1}(\tilde{M}_2^N, \hat{M}_2^N) \Rightarrow 0 \qquad \text{as } N \to \infty,$

then by Theorem 11.4.7 in [20] the proof will be complete. However, (5.14) follows by (5.13), Proposition 5.1 and Theorem 11.4.8 in [20]. The proof is now complete. □

We are now ready to state the main result of this section. Assume first that $\tilde{Q}_0^N \Rightarrow \tilde{Q}_0$ as $N \to \infty$ and let

(5.15) $\qquad \tilde{H} = \tilde{Q}_0 \bar{F}_e \quad \text{and} \quad \tilde{M}_Q = \tilde{Q}_0^+ (G - \bar{F}_e).$

Next, set

(5.16) $\qquad \tilde{M}_1(t) = \int_0^t G(t-s) \, d\tilde{\xi}(s), \qquad t \geq 0,$

and $\tilde{M}_1 = \{\tilde{M}_1(t), t \geq 0\}$, where the process $\tilde{\xi}$ appearing in (5.16) is the limiting process appearing in (5.2) at the beginning of this section. Note also that the integral above may be interpreted as the result of integration by parts.

Next, let $\beta F_e$ be the process $\{\beta F_e(t), t \geq 0\}$, where we recall the definition of $F_e$ from (5.4) above.

Finally, let

(5.17) $\qquad \tilde{Q}_I = \tilde{H} + \tilde{W}_0(F_e) + \tilde{M}_1 + \tilde{M}_2.$

Note that by Theorem 3 of [14], $\tilde{Q}_I$ is the limiting queue length process associated with a sequence of $G/GI/\infty$ queues with identical arrival processes and service time sequence as our original sequence of $G/GI/N$ queues and with $Q_0^N$ customers in service at time zero with residual service time distribution $F_e$.

The following is now our second main result.

THEOREM 5.1. *If the residual service time distribution $F_0 = F_e$ and $\tilde{Q}_0^N \Rightarrow \tilde{Q}_0$ as $N \to \infty$, then $\tilde{Q}^N \Rightarrow \varphi_F^0(\tilde{M}_Q + \tilde{Q}_I - \beta F_e)$ as $N \to \infty$.*

PROOF. Let $f : \mathbb{R} \times D^3 \mapsto \mathbb{R} \times D^3$ be the map defined for $(x_1, x_2, x_3, x_4) \in \mathbb{R} \times D^3$ by

(5.18) $$f((x_1, x_2, x_3, x_4)) = (f_1(x_1), f_2(x_2), f_3(x_3), f_4(x_4)),$$

where $f_1(x_1) = x_1, f_2(x_2) = x_2, f_4(x_4) = x_4$ and

(5.19) $$f_3(x_3)(\cdot) = \int_0^\cdot G(\cdot - s) \, dx_3(s),$$

where the above integral above may be interpreted as the result of integration by parts.

Next, note that since by assumption $F_0 = F_e$ and $\tilde{Q}_0^N \Rightarrow \tilde{Q}_0$ as $N \to \infty$, it follows by Proposition 5.2 that

(5.20) $$(\tilde{Q}_0^N, \tilde{W}_0^N, \tilde{A}^N, \tilde{M}_2^N) \Rightarrow (\tilde{Q}_0, \tilde{W}_0(F_e), \tilde{\xi}, \tilde{M}_2)$$

in $(\mathbb{R} \times D^3[0,\infty), |\cdot| \times d_{J_1}^3)$ as $N \to \infty$, where each of the limiting processes above are independent of one another. Furthermore, we have that each of the limiting processes above are $\mathbb{P}$-a.s. continuous. Thus, since by Lemma A.9 of the Appendix, $f : (\mathbb{R} \times D^3[0,\infty), |\cdot| \times d_{J_1}^3) \mapsto (\mathbb{R} \times D^3[0,\infty), |\cdot| \times d_{J_1}^3)$ is continuous at continuous limit points $(x_1, x_2, x_3, x_4) \in \mathbb{R} \times D^3[0,\infty)$ such that $x_2, x_3, x_4 \in C[0,\infty)$, we have that

(5.21) $$P((\tilde{Q}_0, \tilde{W}_0(F_e), \tilde{\xi}, \tilde{M}_2) \in Disc(f)) = 0.$$

Now note by (5.1), (5.7), (5.9) and the definition of $f$ in (5.18) and (5.19), we have the representation

$$(\tilde{Q}_0^N, \tilde{W}_0^N, \tilde{M}_1^N, \tilde{M}_2^N) = f((\tilde{Q}_0^N, \tilde{W}_0^N, \tilde{A}^N, \tilde{M}_2^N)).$$

It therefore follows by (5.20), the measurability of $f : (\mathbb{R} \times D^3[0,\infty), \mathcal{B}(\mathbb{R}) \times \mathcal{D}^3) \mapsto (\mathbb{R} \times D^3[0,\infty), \mathcal{B}(\mathbb{R}) \times \mathcal{D}^3)$ by Lemma A.9 in the Appendix, (5.21), the Continuous Mapping theorem [1] and the definition $\tilde{M}_1$ in (5.16) that

$$(\tilde{Q}_0^N, \tilde{W}_0^N, \tilde{M}_1^N, \tilde{M}_2^N) \Rightarrow (\tilde{Q}_0, \tilde{W}_0(F_e), \tilde{M}_1, \tilde{M}_2)$$

in $(\mathbb{R} \times D^3[0,\infty), \mathcal{B}(\mathbb{R}) \times \mathcal{D}^3)$ as $N \to \infty$, where each of the limiting processes appearing on the right-hand side above are independent of one another.

Next, since $(\mathbb{R}, \mathbb{R})$ and $(D, \mathcal{D})$ are both separable spaces, it follows by Theorem 11.4.1 in [20] that $\mathbb{R} \times D^3$ is separable under the product topology induced by the maximum metric $|\cdot| \times d_{J_1}^3$. Thus, by the Skorohod representation theorem [20], there exists some alternate probability space, $(\hat{\Omega}, \hat{\mathcal{F}}, \hat{P})$, on which are defined a sequence of processes

$$\{(\hat{Q}_0^N, \hat{W}_0^N, \hat{M}_1^N, \hat{M}_2^N), N \geq 1\},$$

where

(5.22) $\quad (\hat{Q}_0^N, \hat{W}_0^N, \hat{M}_1^N, \hat{M}_2^N) \stackrel{d}{=} (\tilde{Q}_0^N, \tilde{W}_0^N, \tilde{M}_1^N, \tilde{M}_2^N), \qquad N \geq 1,$

and also processes

(5.23) $\quad (\hat{Q}_0, \hat{W}_0(F_e), \hat{M}_1, \hat{M}_2) \stackrel{d}{=} (\tilde{Q}_0, \tilde{W}_0(F_e), \tilde{M}_1, \tilde{M}_2),$

such that

(5.24) $\quad (\hat{Q}_0^N, \hat{W}_0^N, \hat{M}_1^N, \hat{M}_2^N) \to (\hat{Q}_0, \hat{W}_0(F_e), \hat{M}_1, \hat{M}_2)$

in $(\mathbb{R} \times D^3[0, \infty), |\cdot| \times d_{J_1}^3)$ $\hat{\mathbb{P}}$-a.s. as $N \to \infty$. Furthermore, since each of the processes appearing on the right-hand side of (5.24) is continuous, we may assume that above convergence also occurs in $(\mathbb{R} \times D^3[0, \infty), |\cdot| \times u^3)$.

Now, set

$$\hat{M}_Q^N = \hat{Q}_0^{N,+}(G - \bar{F}_e)$$

and

$$\hat{M}_Q = \hat{Q}_0^+(G - \bar{F}_e).$$

It is then clear that

$$\sup_{0 \leq t \leq T} |\hat{M}_Q^N(t) - \hat{M}_Q(t)| \leq |\hat{Q}_0^N - \hat{Q}_0| \sup_{0 \leq t \leq T} |G(t) - \bar{F}_e(t)|$$

$$\leq 2|\hat{Q}_0^N - \hat{Q}_0|,$$

and so it follows by (5.24) that

(5.25) $\quad \hat{M}_Q^N \to \hat{M}_Q \quad$ in $(D[0, \infty), u)$ $\hat{\mathbb{P}}$-a.s. as $N \to \infty$.

Next, letting

$$\hat{H}^N = \hat{Q}_0^N \bar{F}_e,$$

a similar augment shows that

(5.26) $\quad \hat{H}^N \to \hat{H} \quad$ in $(D[0, \infty), u)$ $\hat{\mathbb{P}}$-a.s. as $N \to \infty,$

where

$$\hat{H} = \hat{Q}_0 \bar{F}_e.$$

Finally, it is clear by the Halfin–Whitt condition (5.3) and the boundedness of $F_e$, that

$$(5.27) \qquad \sqrt{N}(1-\rho^N)F_e \to \beta F_e \qquad \text{in } (D[0,\infty), u) \ \hat{\mathbb{P}}\text{-a.s.}$$

as $N \to \infty$.

Thus, letting

$$\hat{Q}_I^N = \hat{H}^N + \hat{W}_0^N + \hat{M}_1^N + \hat{M}_2^N,$$

we have by (5.24)–(5.27) that

$$(5.28) \qquad \hat{M}_Q^N + \hat{Q}_I^N - \sqrt{N}(1-\rho^N)F_e \to \hat{M}_Q + \hat{Q}_I - \beta F_e$$

in $(D[0,\infty), u)$ $\hat{\mathbb{P}}$-a.s. as $N \to \infty$, where

$$\hat{Q}_I = \hat{H} + \hat{W}(F_e) + \hat{M}_1 + \hat{M}_2.$$

Furthermore, it follows by (5.23) that

$$(5.29) \quad \hat{M}_Q^N + \hat{Q}_I^N - \sqrt{N}(1-\rho^N)F_e \stackrel{d}{=} \tilde{M}_Q^N + \tilde{Q}_I^N - \sqrt{N}(1-\rho^N)F_e$$

for $N \geq 1$.

Now set

$$(5.30) \qquad \hat{Q}^N = \varphi_F^0(\hat{M}_Q^N + \hat{Q}_I^N - \sqrt{N}(1-\rho^N)F_e).$$

Since by Proposition 3.1, the map $\varphi_F^0 : (D[0,\infty), \mathcal{D}) \mapsto (D[0,\infty), \mathcal{D})$ is measurable, it follows by (5.11), (5.29) and (5.30) that

$$(5.31) \qquad \hat{Q}^N \stackrel{d}{=} \tilde{Q}^N, \qquad N \geq 1.$$

Furthermore, by the continuity portion of Proposition 3.1 and (5.28),

$$(5.32) \ \hat{Q}^N = \varphi_F^0(\hat{M}_Q^N + \hat{Q}_I^N - \sqrt{N}(1-\rho^N)F_e) \to \varphi_F^0(\hat{M}_Q + \hat{Q}_I - \beta F_e)$$

in $(D[0,\infty), u)$ $\hat{\mathbb{P}}$-a.s. as $N \to \infty$. Since convergence in $(D[0,\infty), u)$ implies convergence in $(D[0,\infty), d_{J_1})$ and almost sure convergence implies convergence in distribution, it now follows by (5.31), the measurability of $\varphi_F^0 : (D[0,\infty), \mathcal{D}) \mapsto (D[0,\infty), \mathcal{D})$ from Proposition 3.1 and (5.32) that $\tilde{Q}^N \Rightarrow \varphi_F^0(\tilde{M}_Q + \tilde{Q}_I - \beta F_e)$ as $N \to \infty$, which completes the proof. □

Note that the diffusion limit for the queue length process given by Theorem 5.1 may be written out in expanded form as the solution to the stochastic convolution equation

$$(5.33) \quad \tilde{Q}_F(t) = \tilde{M}_Q(t) + \tilde{Q}_I(t) - \beta F_e(t) + \int_0^t \tilde{Q}_F^+(t-s) \, dF(s), \qquad t \geq 0.$$

In this representation, we see that $\tilde{Q}_F$ may be decomposed into four separate parts. The second term on the right-hand side of (5.33), $\tilde{Q}_I$, defined by (5.17) above, is the diffusion limit for the $G/GI/\infty$ queue with the same number of initial customers as in our $G/GI/N$ queue and with the same sequence of arrival processes and identical service time distribution as in our $G/GI/N$ queue. It is the primary stochastic component which drives the convolution equation above. The third term on the right-hand side above, $\beta F_e(t)$, arises out of the Halfin–Whitt condition (5.3). The first term on the right-hand side above, $\tilde{M}_Q$, takes into account the discrepancy between the $G/GI/N$ and $G/GI/\infty$ queue in the initial number of customers in the system at time $0-$ who remain in the system at time $t$. Finally, similar to as in the fluid limit of Section 4, the integral term on the right-hand side of (5.33) represents an adjustment term to the infinite server queue limit, $\tilde{Q}_I$ which takes into account the waiting times of the customers in the $G/GI/N$ queue. Note also that the adjustment integral term is positive as one would expect since the number of customers in the $G/GI/N$ queue will always be stochastically larger than in a corresponding $G/GI/\infty$ queue. Furthermore, only the positive portion of the queue length process, $\tilde{Q}_F^+$, appears in the limit since it is only when there are more customers than servers in the system that the finite server approximation to the infinite server queue will be off.

5.4. *The virtual waiting time process.* In this subsection, we study the diffusion scaled virtual waiting time process and the diffusion scaled customer waiting time process for the $G/GI/N$ queue in the Haflin–Whitt regime. For each $t \geq 0$, let $V^N(t)$ denote the hypothetical amount of time that a customer arriving to the $N$th system at time $t$ would have to wait before being served and, for $i \geq 1$, let $V_i^N$ denote the waiting time of the $i$th customer to arrive the system after time $0-$. Note that denoting by $D^N(t)$ the number of departures from the $N$th system by time $t \geq 0$, we have that

$$V^N(t) + t = \inf\{s \geq 0 : D^N(s) \geq A^N(t) + Q_0^N\}$$

for $t \geq 0$ and $V_i^N = V^N(\tau_i^N -)$. In this subsection, we obtain asymptotic results for the diffusion scaled virtual waiting time process $\tilde{V}^N(t) = \{N^{1/2}V^N(t), t \geq 0\}$ and the also for the diffusion scaled customer waiting time process $\tilde{\hat{V}}^N(t) = \{N^{1/2}V_{\lfloor Nt \rfloor}^N, t \geq 0\}$. Our main approach will be to leverage off of the results of Section 5.3 on the diffusion scaled queue length process.

The main result of this subsection is the following theorem which provides a weak limit for the sequences of diffusion scaled virtual waiting time and customer waiting time processes in the Halfin–Whitt regime. Note that as a byproduct, this limit implies that waiting times in the Halfin–Whitt regime are of order $N^{-1/2}$ and thus decrease as the number of servers become large.

Furthermore, we also note that as expounded upon further below, the limit we obtain is similar in form to that informally obtained by Mandelbaum and Momčilović in [16] for the diffusion scaled virtual waiting time process. We then have the following.

PROPOSITION 5.3. *If $A^N(0) = 0$ $\mathbb{P}$-a.s. for each $N \geq 1$, and the residual service time distribution $F_0 = F_e$ and $\tilde{Q}_0^N \Rightarrow \tilde{Q}_0$ as $N \to \infty$, then $\tilde{V}^N \Rightarrow \tilde{V}$ as $N \to \infty$ and $\tilde{\tilde{V}}^N \Rightarrow \tilde{V}$ as $N \to \infty$, where $\tilde{V}$ is given by the unique solution to the integral equation*

$$(5.34) \quad \tilde{V}(t) = \left( \tilde{M}_Q(t) + \tilde{Q}_I(t) - \beta F_e(t) + \int_0^t \tilde{V}(t-s) \, dF(s) \right)^+, \qquad t \geq 0.$$

In [16], in the case of renewal arrivals, the diffusion scaled limiting stationary virtual waiting time process was conjectured to be expressed in terms of the supremum over an infinite weighted full $K$-ary tree

$$(5.35) \qquad \tilde{V}_\infty(t) = \sup_{T \subset F[t]} (W_T)^+, \qquad t \in \mathbb{R}.$$

Furthermore, it was shown in Lemma 10 of [16] that the process $\tilde{V}_\infty = \{\tilde{V}_\infty(t), t \geq 0\}$ defined by (5.35) satisfies the stochastic integral equation

$$(5.36) \quad \tilde{V}_\infty(t) = \left( \tilde{X}(t) - \beta + \int_0^\infty \tilde{V}(t-u) \, dF(u) \right)^+, \qquad t \in \mathbb{R},$$

where $\tilde{X} = \{\tilde{X}(t), t \geq 0\}$ is a stationary Gaussian process whose covariance function may be explicitly calculated.

We now show informally by taking the limit as $t$ approaches $\infty$ on both sides of (5.34) that (5.36) may actually be viewed as the stationary version of the limiting diffusion scaled virtual waiting time process. First, note that by the definition of $\tilde{M}_Q$ in (5.15), we have that $\tilde{M}_Q(t) \to 0$ $\mathbb{P}$-a.s. as $t \to \infty$. Next, since $F_e$ is a distribution function, it follows that $\beta F_e(t) \to \beta$ $\mathbb{P}$-a.s. as $t \to \infty$. Thus, for $t$ large, we have from (5.34) that

$$(5.37) \quad \tilde{V}(t+s) \approx \left( \tilde{Q}_I(t+s) - \beta + \int_0^{t+s} \tilde{V}(t+s-u) \, dF(u) \right)^+$$

for $s \geq 0$, where the $\approx$ sign is meant to mean approximately equal. One may next check that for $t$ large the process $\{\tilde{Q}_I(t+s), s \geq 0\}$ is approximately equal in law to that of $\tilde{X}$ in (5.37). Finally, for $t$ large, the integral on the right-hand side of (5.37) may be taken to be an interval over the entire nonnegative portion of the real line. Thus, setting $\tilde{V}_t = \{\tilde{V}(t+s), s \geq 0\}$, we obtain from (5.37) and the preceding discussion that

$$(5.38) \quad \tilde{V}_t(s) \approx \left( \tilde{X}(s) - \beta + \int_0^\infty \tilde{V}_t(s-u) \, dF(u) \right)^+, \qquad s \geq 0,$$

which is the similar to (5.38). As $t$ approaches $\infty$, the approximation in (5.38) becomes more precise and indeed a completely rigorous argument of the discussion above may be given but we omit the details here.

In order to prove Proposition 5.3, we must first provide an intermediary result for the diffusion scaled number of customers waiting to be served. Note that since the system is operating under a nonidling policy, the total number of customers waiting to be served at time $t$ is given by the quantity $(Q^N(t) - N)^+$. We also define

$$\tilde{Q}^{N,+}(t) = \frac{(Q^N(t) - N)^+}{N^{1/2}}, \qquad t \geq 0,$$

to be the diffusion scaled number of customers waiting to be served at time $t$ and set $\tilde{Q}^{N,+} = \{N^{-1/2}(Q^N(t) - N)^+, t \geq 0\}$ to be the diffusion scaled number of customers waiting to be served process. Note also the relationship $\tilde{Q}^{N,+}(t) = (\tilde{Q}^N(t))^+$ which will be taken advantage of in the following result.

PROPOSITION 5.4. *If the residual service time distribution $F_0 = F_e$ and $\tilde{Q}_0^N \Rightarrow \tilde{Q}_0$ as $N \to \infty$, then we have the joint convergence $(\tilde{Q}^N, \tilde{Q}^{N,+}) \Rightarrow (\tilde{Q}_F, \tilde{Q}_F^+)$ in $(D^2[0,\infty), \mathcal{D}^2)$ as $N \to \infty$.*

PROOF. Define the function $f:(D[0,\infty), d_{J_1}) \mapsto (D^2[0,\infty), d_{J_1}^2)$ by $f(x_1) = (x_1, x_1^+)$ for $x_1 \in D[0,\infty)$, where $x_1^+ = (\max\{0, x_1(t)\}, t \geq 0)$. We now claim that the function $f:(D[0,\infty), d_{J_1}) \mapsto (D^2[0,\infty), d_{J_1}^2)$ is continuous. Assume that $x_1^n \to x_1$ in $d_{J_1}$. We then have that

$$\begin{aligned} d_{J_1}^2(f(x_1^n), f(x_1)) &= d_{J_1}^2((x_1^n, x_1^{n,+}), (x_1, x_1^+)) \\ &= \max\{d_{J_1}(x_1^n, x_1), d_{J_1}(x_1^{n,+}, x_1^+)\} \\ &\leq d_{J_1}(x_1^n, x_1) \\ &\to 0 \qquad \text{as } n \to \infty, \end{aligned}$$

and thus the claim is proven. The proof now follows by the representation $(\tilde{Q}^N, \tilde{Q}^{N,+}) = f(\tilde{Q}^N)$, the Continuous Mapping theorem [20] and Theorem 5.1 above. □

Using Lemma A.10 in the Appendix, we may now prove the proposition above.

PROOF OF PROPOSITION 5.3. By Proposition 5.4, $\tilde{Q}^{N,+} \Rightarrow \tilde{Q}^+$ as $N \to \infty$ and by assumption we have that $\tilde{A}^N \Rightarrow \tilde{\xi}$ as $N \to \infty$. Furthermore, by (5.1) and the heavy traffic condition (5.3), $\rho^N \to 1$ as $N \to \infty$. Thus, by Lemma A.10 in the Appendix, $\tilde{V}^N \Rightarrow \tilde{V} = \tilde{Q}^+$ as $N \to \infty$, which by the representation (5.33), completes the proof. □

5.5. *An alternative representation.* The representation of the limiting diffusion scaled queue length process given by (5.33) provides a convenient expression for $\tilde{Q}_F$ in terms of a corresponding infinite server queue limit, $\tilde{Q}_I$. However, it is not evident that the limiting process of (5.33) is equivalent to the diffusion limit of Theorem 2 of Halfin and Whitt [8] for the specific case of exponentially distributed service times. The following corollary provides an alternative representation of the limit of (5.33) which may be used to rigorously verify this equivalency. Let $M = \{M(t), t \geq 0\}$ be the renewal function associated with the pure renewal process with interarrival distribution given by the service time distribution $F$. Recall that for $t \geq 0$, $M(t)$ is by definition equal to the expected number of renewals by time $t$. Furthermore, by Exercise 3.4 of [19], $M$ is given by the unique solution to the renewal equation

$$(5.39) \qquad M(t) = F(t) + \int_0^t M(t-u) \, dF(u) \qquad \text{for } t \geq 0.$$

Let us also set

$$(5.40) \qquad \tilde{\zeta}(t) = \tilde{M}_Q(t) + \tilde{Q}_I(t), \qquad t \geq 0,$$

and $\tilde{\zeta} = \{\tilde{\zeta}(t), t \geq 0\}$. We then have the following result.

COROLLARY 5.2. *The limiting process, $\tilde{Q}_F$, of Theorem 5.1 may be equivalently expressed as the unique strong solution to*

$$(5.41) \quad \tilde{Q}_M(t) = \tilde{\zeta}(t) + \int_0^t \tilde{\zeta}(t-u) \, dM(u) - \beta t - \int_0^t \tilde{Q}_M^-(t-u) \, dM(u)$$

*for $t \geq 0$, where $\tilde{Q}_M^-(t) = \min(\tilde{Q}_M(t), 0)$.*

PROOF. Let $F = \{F(t), t \geq 0\}$ be a distribution function and $r = \{r(t), t \geq 0\}$ be an unknown function satisfying the integral equation of renewal type,

$$(5.42) \qquad r(t) = H(t) + \int_0^t r(t) \, dF(t-u) \qquad \text{for } t \geq 0$$

for some $H = \{H(t), t \geq 0\}$. If $H$ is a locally bounded function, then (5.42) has a unique locally bounded solution [11], which is given by

$$(5.43) \qquad r(t) = H(t) + \int_0^t H(t-u) \, dM(u),$$

where $M = \{M(t), t \geq 0\}$ is the solution to the renewal equation (5.39).

By the definition of $\tilde{\zeta}$ in (5.40) and the representation (5.33), the limiting process of Theorem 5.1 may be written as

$$\tilde{Q}_F(t) = \tilde{\zeta}(t) - \beta F_e(t) + \int_0^t \tilde{Q}_F^+(s) \, dF(t-s) \qquad \text{for } t \geq 0.$$

Next, since $\tilde{Q}_F = \tilde{Q}_F^+ + \tilde{Q}_F^-$, we have

$$\tilde{Q}_F^+(t) = \tilde{\zeta}(t) - \beta F_e(t) - \tilde{Q}_F(t)^- - \int_0^t \tilde{Q}_F^+(s)\,dF(t-s).$$

Furthermore, it follows that $\tilde{\zeta} - \beta F_e + \tilde{Q}_F^-$ is almost surely a locally bounded function since it is almost surely an element of $D[0, \infty)$. It therefore follows from (5.42) and (5.43) that

$$\tilde{Q}_F^+(t) = \tilde{\zeta}(t) - \beta F_e(t) - \tilde{Q}_F^-(t) + \int_0^t \tilde{\zeta}(t-u)\,dM(u)$$
$$- \beta \int_0^t F_e(t-u)\,dM(u) - \int_0^t \tilde{Q}_F^-(t-u)\,dM(u),$$

or, equivalently,

$$(5.44) \quad \begin{aligned} \tilde{Q}_F(t) &= \tilde{\zeta}(t) + \int_0^t \tilde{\zeta}(t-u)\,dM(u) \\ &\quad - \beta\left(F_e(t) + \int_0^t F_e(t-s)\,dM(s)\right) \\ &\quad - \int_0^t \tilde{Q}_F^-(t-u)\,dM(u). \end{aligned}$$

However, since

$$F_e(t) + \int_0^t F_e(t-s)\,dM(s) = t, \qquad t \geq 0,$$

if follows by (5.44) that

$$(5.45) \quad \tilde{Q}_F(t) = \tilde{\zeta}(t) + \int_0^t \tilde{\zeta}(t-u)\,dM(u) - \beta t - \int_0^t \tilde{Q}_F^-(t-u)\,dM(u),$$

which completes the proof. □

Note that Corollary 5.2 also implies the convergence $\tilde{Q}^N \Rightarrow \tilde{Q}_M$ as $N \to \infty$, where $\tilde{Q}_M$ is the stochastic process given by the unique strong solution to (5.41). The proof of this result is trivial and essentially proceeds in two stages. First, one may use Theorem 5.1 to show the weak convergence $\tilde{Q}^N \Rightarrow \tilde{Q}_F$ as $N \to \infty$ and then use Corollary 2 to show the equivalency in law between $\tilde{Q}_F$ and $\tilde{Q}_M$. In a sequel [18] to this paper, we provide a more direct proof of the convergence of $\tilde{Q}^N$ to $\tilde{Q}_M$. This follows along more traditional line of conventional heavy-traffic proofs in where the queue length process is modeled by a conservation of flow as the difference between the number of arrivals and the number of departures. The integral term on the right-hand side of (5.41) then turns out to representative of the idle time processes of the servers in the system.

As noted above, in the case of the $GI/M/N$ queue, Corollary 5.2 may be used to obtain the original diffusion limit result provided by Halfin and Whitt [6]. This may be seen by first noting that for exponentially distributed service times, the renewal function $M$ in (5.41) is the renewal function for a rate 1 Poisson process, which is simply given by $M(t) = t$. Thus, the limiting process of Corollary 5.2 may be written as

$$(5.46) \quad \tilde{Q}_M(t) = \tilde{\zeta}(t) + \int_0^t \tilde{\zeta}(s)\,ds - \beta t - \int_0^t \tilde{Q}_M^-(s)\,ds, \qquad t \geq 0.$$

Furthermore, using (5.40), (5.15) and (5.17), extensive covariance calculations show that

$$B(t) = \tilde{\zeta}(t) + \int_0^t \tilde{\zeta}(s)\,ds, \qquad t \geq 0,$$

is a Brownian motion with infinitesimal variance $1 + \sigma^2$, where $\sigma^2$ is the variance of the interarrival times. Therefore, the process (5.46) is a diffusion with infinitesimal drift $m(x) = -\beta$ for $x \geq 0$ and $m(x) = -x - \beta$ for $x < 0$ and infinitesimal variance $1 + \sigma^2$, which is in agreement with Theorem 3 of Halfin and Whitt [6].

**6. Conclusion.** In this paper, we have studied the $G/GI/N$ queue in the Halfin–Whitt regime. In our first main result, we obtained a first-order approximation to the queue length process. This approximation turned out to be the solution to a nonlinear convolution equation. Next, after centering the queue length process by its deterministic fluid limit and scaling by an appropriate constant, we obtained a second-order stochastic approximation as well. Our limiting stochastic process is nonlinear, stochastic convolution equation which is driven by a Gaussian process and includes a drift term which incorporates a time lag. In the case of exponentially distributed service times, it can be shown that this process is equivalent to the diffusion process obtained by Halfin and Whitt [6].

In the sequel to this paper [18], we provide a direct approach to the proof of Corollary 5.2. This is based off of a conservation of flow equation where we write the queue length process as the difference between the number of arrivals and the number of departures. In this case, central limit theorem type results for sums of i.i.d. renewal process will have to be invoked.

In the future, it would be nice to have a better understanding of the limiting process we have obtained. Ideally, one would like to solve for its limiting distribution. Unfortunately, this in general appears to be a difficult problem. Therefore, if analytical solutions cannot be found, efficient numerical procedures might perhaps be developed. Simulation studies could also be conducted to test the accuracy of the proposed approximations relative to their actual values. This would be especially interesting when the system is close to being in the Halfin–Whitt regime.

## APPENDIX

In the appendix, we provide the proofs of Propositions 3.1, 4.2 and 5.1 and Lemmas A.9 and A.10. We begin with the proof of Proposition 3.1.

PROOF OF PROPOSITION 3.1. Suppose first that $B$ is concentrated on the point $c > 0$. In this case, it is clear that the solution to (3.1) satisfies the recursion

$$z(t) = x(t), \qquad 0 \leq t < c, \tag{A.1}$$

and

$$z(t) = x(t) + (z(t-c) + a)^+, \qquad t \geq c,$$

in which case it is clearly unique. Furthermore, defining $\varphi_B^a : D[0,\infty) \mapsto D[0,\infty)$ to be the solution to this recursion, it follows that

$$\|\varphi_B^a(x_1) - \varphi_B^a(x_2)\|_t = \|x_1 - x_2\|_t$$

for $0 \leq t < c$. Now suppose that for some integer $k$, we have

$$\|\varphi_B^a(x_1) - \varphi_B^a(x_2)\|_t \leq k\|x_1 - x_2\|_t \tag{A.2}$$

for $(k-1)c \leq t < kc$. It then follows that for $k < t \leq (k+1)c$,

$$\begin{aligned}\|\varphi_B^a(x_1) - \varphi_B^a(x_2)\|_t &\leq \|x_1 - x_2\|_t + \|\varphi_B^a(x_1) - \varphi_B^a(x_2)\|_{t-c} \\ &\leq \|x_1 - x_2\|_t + k\|x_1 - x_2\|_t \\ &= (k+1)\|x_1 - x_2\|_t.\end{aligned}$$

By induction, this implies that the relationship (A.2) must hold for all $t$, which show that $\varphi_B^a$ is Lipschitz continuous if $B$ is concentrated on a single point. The proof of measurability of $\varphi_B^a$ for the case of $B$ concentrated at a single point will be included below.

Suppose now that there exists a $\delta > 0$ such that $B(y+\delta) - B(y) < \varepsilon$ for some $0 < \varepsilon < 1$ for all $y \geq 0$. Such a $\delta$ will always exist so long as $B$ is not concentrated on a single point. We now provide proofs of existence, uniqueness and Lipschitz continuity for this case.

*Existence:* We use the method of successive approximations. Let $u_0 = 0$ and recursively define

$$u_{n+1}(t) = x(t) + \int_0^t (u_n(t-s) + a)^+ \, dB(s), \qquad t \geq 0, \tag{A.3}$$

for $n \geq 1$ and note that

$$\begin{aligned}u_{n+1}(t) &- u_n(t) \\ &= \int_0^t ((u_n(t-s) + a)^+ - (u_{n-1}(t-s) + a)^+) \, dB(s), \qquad t \geq 0.\end{aligned} \tag{A.4}$$

We now show by induction that for each integer $1 \leq k \leq \lceil \delta^{-1} T \rceil$,

(A.5) $\qquad \|u_{n+1} - u_n\|_{j\delta} < j^j n^j \varepsilon^n \|x\|_T \qquad$ for $j = 1, \ldots, k,$

for $n \geq 1$. For the base case $k = 1$, observe by (A.4) that

$$\|u_{n+1} - u_n\|_\delta \leq B(\delta) \|u_n - u_{n-1}\|_\delta < \varepsilon \|u_n - u_{n-1}\|_\delta,$$

and so, since

(A.6) $\qquad \|u_1 - u_0\|_\delta = \|x_1 - 0\|_\delta \leq \|x_1\|_T,$

we have the relationship

(A.7) $\qquad \|u_{n+1} - u_n\|_\delta < \varepsilon^n \|x\|_T \leq n\varepsilon^n \|x\|_T$

for $n \geq 1$. We now proceed with the induction step. Assume that (A.5) holds for some $k$ and we will now show that it holds for $(k+1)$. Note that by (A.4) and the induction hypothesis (A.5),

(A.8) $\qquad \|u_{n+1} - u_n\|_{(k+1)\delta}$

$$\leq \sum_{j=1}^k \varepsilon \|u_n - u_{n-1}\|_{j\delta} + \varepsilon \|u_n - u_{n-1}\|_{(k+1)\delta}$$

$$\leq \sum_{j=1}^k \varepsilon j^j (n-1)^j \varepsilon^{n-1} \varepsilon \|x\|_T + \varepsilon \|u_n - u_{n-1}\|_{(k+1)\delta}$$

(A.9) $\qquad \leq k^{k+1} n^k \varepsilon^n \|x\|_T + \varepsilon \|u_n - u_{n-1}\|_{(k+1)\delta}.$

Furthermore, since as in (A.6),

$$\|u_1 - u_0\|_{l\delta} \leq \|x\|_T$$

for all $l = 1, \ldots, \lceil \delta^{-1} T \rceil$, it follows by repeated iteration of (A.9) that

$$\|u_{n+1} - u_n\|_{(k+1)\delta} \leq k^{k+1} \left( \sum_{i=0}^n i^k \right) \varepsilon^n \|x\|_T$$

$$\leq k^{k+1} (n^{k+1} + 1) \varepsilon^n \|x\|_T$$

$$\leq (k+1)^{k+1} n^{k+1} \varepsilon^n \|x\|_T,$$

and so the induction hypothesis has been proved.

Thus, since

$$\sum_{n=1}^\infty \|u_{n+1} - u_n\|_T \leq \sum_{n=1}^\infty \|u_{n+1} - u_n\|_{\lceil \delta^{-1} T \rceil \delta}$$

$$\leq \sum_{n=1}^\infty \lceil \delta^{-1} T \rceil^{\lceil \delta^{-1} T \rceil} n^{\lceil \delta^{-1} T \rceil} \varepsilon^n \|x\|_T$$

$$< \infty,$$

it follows that $\{u_n\}$ is a Cauchy sequence. Furthermore, since the space $D[0,\infty)$ is a Banach space under the supremum metric $u$, there exists a limit point $u^\star$ of $\{u_n\}$. Taking limits on both sides of (A.3), we now see that $u^\star$ is a solution to (3.1), which completes the proof of existence for the case of nondegenerate distributions.

*Uniqueness:* Suppose that $u$ and $v$ both satisfy (3.1) and let

$$\Delta(t) = u(t) - v(t) = \int_0^t ((u(t-s)+a)^+ - (v(t-s)+a)^+) \, dB(s), \qquad t \geq 0.$$

We then have for $0 \leq t \leq \delta$, that

$$|\Delta(t)| \leq \int_0^t |(u(t-s)+a)^+ - (v(t-s)+a)^+| \, dB(s) \leq \varepsilon \|\Delta\|_\delta,$$

which implies that $\Delta(t) = 0$ on $[0, \delta]$. Next, for $\delta < t \leq 2\delta$, we have that

$$|\Delta(t)| \leq \varepsilon \|\Delta\|_\delta + \varepsilon \|\Delta\|_{2\delta} = \varepsilon \|\Delta\|_{2\delta},$$

which implies that $\Delta(t) = 0$ on $[\delta, 2\delta]$. Iterating the above argument until we reach $T$ completes the proof.

*Lipschitz continuity:* Note that for $0 \leq t < \delta$, we have

$$\|\varphi_B^a(x_2) - \varphi_B^a(x_1)\|_\delta \leq \|x_2 - x_1\|_\delta + \varepsilon \|\varphi_B^a(x_2) - \varphi_B^a(x_1)\|_\delta,$$

which implies that

$$\|\varphi_B^a(x_2) - \varphi_B^a(x_1)\|_\delta \leq (1-\varepsilon)^{-1} \|x_2 - x_1\|_\delta.$$

Next, for $\delta < t \leq 2\delta$, we have

$$\|\varphi_B^a(x_2) - \varphi_B^a(x_1)\|_{2\delta}$$
$$\leq \|x_2 - x_1\|_{2\delta} + \varepsilon \|\varphi_B^a(x_2) - \varphi_B^a(x_1)\|_\delta + \varepsilon \|\varphi_B^a(x_2) - \varphi_B^a(x_1)\|_{2\delta}$$
$$\leq \frac{1}{(1-\varepsilon)} \|x_1 - x_2\|_{2\delta} + \varepsilon \|\varphi_B^a(x_2) - \varphi_B^a(x_1)\|_{2\delta},$$

which implies that

$$\|\varphi_B^a(x_2) - \varphi_B^a(x_1)\|_{2\delta} \leq \frac{1}{(1-\varepsilon)^2} \|x_1 - x_2\|_{2\delta}.$$

Iterating the above argument for $k = \lceil \delta^{-1} T \rceil - 2$ more time intervals completes the proof.

Finally, we provide a proof of measurability of $\varphi_B^a$ for the case of a general $B$.

*Measurability:* We begin by defining the function $\Psi_B^a : D[0,\infty) \to D : [0,\infty)$ by

$$\Psi_B^a(u)(t) = \int_0^t (u(t-s)+a)^+ \, dB(s), \qquad t \geq 0.$$

We now show that $\Psi_B^a$ is measurable with respect to the Borel $\sigma$-field $\mathcal{D}$ generated by the Skorohod $J_1$ topology. Note that since $\mathcal{D}$ is equal to the Kolmogorov $\sigma$-field, which is generated by the finite-dimensional cylinder sets, it is sufficient to check that for each $n \geq 1$ and $A_1, A_2, \ldots, A_n \in B(\mathbb{R})$,

$$\{u \in D[0, \infty) : (\Psi_B^a(u)(t_1), \ldots, \Psi_B^a(u)(t_n)) \in (A_1, \ldots, A_n)\} \in \mathcal{D}$$

for $0 \leq t_1 < t_2 < \cdots < t_n$. However, since $\sigma$-algebras are closed under finite intersections, it is sufficient to check that for each $t \geq 0$, $\Psi_B^a(\cdot)(t)$ is measurable. In order to show this, we first decompose $B$ into its continuous and discrete parts so that

$$B(t) = B_c(t) + B_d(t), \qquad t \geq 0,$$

where we write

(A.10) $$B_d(t) = \sum_{n=1}^{\infty} c_n \delta_{(p_n)}(t)$$

for the discrete part of $B$. We then show that both $\Psi_{B_c}^a$ and $\Psi_{B_d}^a$ are measurable functions and so, since the sum of two measurable functions from $(D, \mathcal{D})$ to $(D, \mathcal{D})$ is measurable, and $\Psi_B^a = \Psi_{B_c}^a + \Psi_{B_d}^a$, we have the desired measurability of $\Psi_B^a$.

We begin with the proof of measurability for $\Psi_{B_c}^a$ for which it will be sufficient to show that for each $t \geq 0$, $\Psi_{B_c}^a(\cdot)(t)$ is continuous when viewed as a function from $(D[0, \infty), d_{J_1})$ to $(\mathbb{R}, |\cdot|)$. Let $u_n \to u$ under the metric $d_{J_1}$. This then implies that $u_n(t) \to u(t)$ for all but a countable number of $t$ (see, for instance, page 247 of [1]). Furthermore, the measure defined by $B_c$ assigns measure 0 to all countable sets. Thus, since for each $t \geq 0$, the sequence $\{\sup_{0 \leq s \leq t} |u_n(s)|, n \geq 1\}$ is bounded, it follows by Theorem 3 of [4] that

(A.11) $$|\Psi_{B_c}^a(u_n)(t) - \Psi_{B_c}^a(u)(t)|$$

$$= \left| \int_0^t ((u_n(t-s) + a)^+ - (u(t-s) + a)^+) \, dB(s) \right|$$

$$\leq \int_0^t |u_n(t-s) - u(s)| \, dB(s)$$

(A.12) $$\to 0 \quad \text{as } n \to \infty.$$

This completes the proof of the measurability of $\Psi_{B_c}^a$.

Now consider $\Psi_{B_d}^a$. It is clear from (A.10) that

$$\Psi_{B_d}^a(u)(t) = \sum_{k=1}^{\infty} \Upsilon_k(u)(t), \qquad t \geq 0,$$

where
$$\Upsilon_k(u)(t) = c_k 1\{t \geq p_k\}(u(t - p_k) + a)^+.$$

For each $n \geq 1$, define
$$\Psi_{B_d}^{a,n}(u)(t) = \sum_{k=1}^{n} \Upsilon_k(u)(t), \qquad t \geq 0.$$

We then have that for each $u \in D[0, \infty)$ and $t \geq 0$,
$$\sup_{0 \leq s \leq t} |\Psi_{B_d^n}^a(u)(s) - \Psi_{B_d}^a(u)(s)|$$
$$= \sup_{0 \leq s \leq t} \left| \sum_{k=n}^{\infty} c_k 1\{s \geq p_k\}(u(s - p_k) + a)^+ \right|$$
$$\leq \sup_{0 \leq s \leq t} |u(s) + a| \sum_{k=n}^{\infty} c_k$$
$$\to 0 \qquad \text{as } n \to \infty,$$

and so it follows that $\Psi_{B_d}^a(u)$ is the pointwise limit in $(D[0, \infty), u)$ of $\Psi_{B_d}^{a,n}(u)$ as $n \to \infty$. Thus, if each $\Psi_{B_d}^{a,n}(u)$ is measurable, it will follow that $\Psi_{B_d}^a$ is measurable as well. However, in order to show that $\Psi_{B_d}^{a,n}$ is measurable, it will suffice to show that each $\Upsilon_k$ is measurable since the sum of a finite number of measurable functions is measurable. The fact that $\Upsilon_k$ is measurable may been seen by noting that $\Upsilon_k$ is first the translation of the function $u$ by a constant $p_k$ and then a multiplication by a constant $c_k$. Both of these functions are easily seen to be measurable functions and so $\Upsilon_k$, being the composition of two measurable functions, is measurable as well. This completes the proof of the measurability of $\Psi_{B_d}^a$.

Now define the map $\Xi_B^a : (D[0, \infty), \mathcal{D}) \mapsto (D[0, \infty), \mathcal{D})$ by
$$\Xi_B^a(u)(t) = x(t) + \Psi_B^a(u)(t), \qquad t \geq 0.$$

It is clear that $\Xi_B^a$ is measurable since $\Psi_B^a$ is measurable. Furthermore, from the existence portion of the arguments above, it follows that for each $x \in D[0, \infty)$,
$$\varphi(x) = \lim_{n \to \infty} \Xi_B^{a,n}(0),$$

where $\Xi_B^{a,n}(x) = \Xi_B^{a,n-1} \circ \Xi_B^a$ is the $n$-fold composition of $\Xi_B^a$ with itself and the limit is taken with respect to the metric of uniform convergence over bounded intervals, $u$. Thus, since the composition of two measurable functions is measurable, it follows that $\Xi_B^{a,n}$ is measurable for each $n$. But this then implies that $\varphi_B^a$, being the pointwise limit of a sequence of measurable functions, is measurable as well, and so the proof is now complete. $\square$

The next portion of the Appendix is devoted to the proofs of Propositions 4.2 and 5.1. Our proof of Proposition 4.2 closely parallels the proofs of Lemmas 3.4–3.8 of [14]. In order to begin, we must first set up the following notation. Let $\hat{A}^N(t)$ be equal to the number of customers in the $N$th system who entered service after time $0-$ but before or at time $t$. We then define the two parameter process

$$(A.13) \quad V^N(t,x) = \sum_{i=1}^{\hat{A}^N(t)} (1\{\eta_i \leq x\} - F(x)), \qquad t \geq 0, x \geq 0,$$

where we recall the definition in Section 2 of $\eta_i$ as the service time of the $i$th customer to arrive to the system after time $0-$. Note that by setting

$$U^N(t,x) = \sum_{i=1}^{\lfloor Nt \rfloor} (1\{F(\eta_i) \leq x\} - x), \qquad t \geq 0, 0 \leq x \leq 1,$$

we have

$$(A.14) \quad V^N(t,x) = U^N(\check{A}^N(t), F(x)),$$

where

$$(A.15) \quad \check{A}^N(t) = \frac{\hat{A}^N(t)}{N}, \qquad t \geq 0.$$

It then follows from the definition of $M_2^N$ in (2.5) that

$$(A.16) \quad M_2^N(t) = \int_0^t \int_0^t 1\{s+x \leq t\} \, dV^N(s,x),$$

where the integrals above are taken over the closed intervals $[0,t]$. We now decompose $M_2^N$ in two processes, $G^N$ and $H^N$. Let

$$(A.17) \quad L^N(t,x) = \sum_{i=1}^{\hat{A}^N(t)} \left( 1\{\eta_i \leq x\} - \int_0^{x \wedge \eta_i} \frac{dF(y)}{1 - F(y-)} \right), \qquad t \geq 0, x \geq 0,$$

where $F(y-) = \lim_{x \to y} F(x)$.

By (A.13) and (A.17), we have that

$$(A.18) \quad V^N(t,x) = -\int_0^x \frac{V^N(t,y-)}{1 - F(y-)} dF(y) + L^N(t,x).$$

Therefore, by (A.16) and (A.18), we have

$$(A.19) \quad M_2^N(t) = G^N(t) + H^N(t),$$

where

$$G^N(t) = -\int_0^t \frac{V^N(t-x, x-)}{1 - F(x-)} dF(x), \qquad t \geq 0,$$

and

$$(A.20) \qquad H^N(t) = \int_0^t \int_0^t 1\{s + x \leq t\} \, dL^N(s, x), \qquad t \geq 0.$$

We set $G^N = \{G^N(t), t \geq 0\}$ and $H^N = \{H^N(t), t \geq 0\}$ and note that (A.19) is the desired decomposition of $M_2^N$. It will be useful in proving several results related to $M_2^N$ such as tightness and weak convergence.

Now, for each $k \geq 1$, let

$$H_k^N(t) = \sum_{i=1}^{\hat{A}^N(t) \wedge k} \left( 1\{0 < \eta_i \leq t - \hat{\tau}_i^N\} - \int_{0+}^{\eta_i \wedge (t - \hat{\tau}_i^N)^+} \frac{dF(u)}{1 - F(u-)} \right)$$

for $t \geq 0$, where

$$\hat{\tau}_i^N = \inf\{t \geq 0 : \hat{A}^N(t) \geq i\}$$

is the time at which the $i$th customer to enter service after time $0-$ begins being served. We also set $H_k^N = \{H_k^N(t), t \geq 0\}$. Furthermore, we define the filtration $\mathbf{H}^N = (\mathcal{H}_t^N, t \geq 0)$ by

$$\mathcal{H}_t^N = \sigma\{Q_0^N\} \vee \sigma\{\tilde{\eta}_i, i \geq 1\} \vee \sigma\{\xi_i, i \geq 1\}$$

$$\vee \sigma\{1\{\eta_i = 0\}, 1\{\eta_i \leq s - \hat{\tau}_i^N\}, s \leq t, i = 1, \ldots, \hat{A}^N(t)\}$$

$$\vee \sigma\{\hat{A}^N(s), s \leq t\} \vee \mathcal{N},$$

where $\xi_i$ is as defined in (2.1) of Section 2 and $\mathcal{N}$ is the $\mathbb{P}$ completion of $\mathcal{F}$. It easy to see that $\mathbf{H}^N$ satisfies the usual conditions and is actually a filtration.

The following lemma is similar to Lemma 3.5 of [14]. Consequently, the proof that follows is a straightforward adaptation of that found in [14].

LEMMA A.1. *The process $H_k^N$ is an $\mathbf{H}^N$-square-integrable martingale with predictable quadratic variation process*

$$\langle H_k^N \rangle(t) = \sum_{i=1}^{\hat{A}^N(t) \wedge k} \int_{0+}^{\eta_i \wedge (t - \hat{\tau}_i^N)^+} \frac{1 - F(u)}{(1 - F(u-))^2} \, dF(u), \qquad t \geq 0.$$

PROOF. We first decompose $H_k^N$ by writing

$$H_k^N(t) = \sum_{i=1}^k H^{N,i}(t),$$

where

$$(A.21) \quad H^{N,i}(t) = 1\{0 < \eta_i \leq t - \hat{\tau}_i^N\} - \int_{0+}^{\eta_i \wedge (t - \hat{\tau}_i^N)^+} \frac{dF(u)}{1 - F(u-)}.$$

As in [14], the proof now proceeds in three parts. They are:

1. For each $i \geq 1$, the process $H^{N,i} = (H^{N,i}(t), t \geq 0)$ is an $\mathbf{H}^N$-square-integrable martingale.
2. The predictable quadratic covariation process of $H^{N,i}$ is given by

$$\langle H^{N,i}\rangle(t) = \int_{0+}^{\eta_i \wedge (t-\hat{\tau}_i^N)^+} \frac{1-F(u)}{(1-F(u-))^2}\, dF(u). \tag{A.22}$$

3. The martingales $H^{N,i}$ and $H^{N,j}$ are orthogonal for $i \neq j$.

These three statements are then sufficient to imply the conclusion of the lemma. We begin with the proof of part 1.

First, note that $H^{N,i}$ is $\mathbf{H}^N$-adapted and, furthermore, we have that

$$\sup_{t \geq 0} E(H^{N,i}(t))^2 < \infty.$$

We now prove the martingale property for $H^{N,i}$ by showing that for $s < t$,

$$1\{\hat{\tau}_i^N > s\} E[H^{N,i}(t)|\mathcal{H}_s^N] = 0 \tag{A.23}$$

and

$$1\{\hat{\tau}_i^N \leq s\} E[H^{N,i}(t)|\mathcal{H}_s^N] = H^{N,i}(s). \tag{A.24}$$

We begin with (A.23). First, note that $\hat{\tau}_i^N$ is an $\mathbf{H}^N$-stopping time since $\sigma(\hat{A}^N(s), s \leq t) \subset \mathcal{H}_t^N$ for each $t \geq 0$ and so we may define the $\sigma$-field $\mathcal{H}_{\hat{\tau}_i^N}^N$. Furthermore,

$$1\{\hat{\tau}_i^N > s\} E[H^{N,i}(t)|\mathcal{H}_s^N] = 1\{\hat{\tau}_i^N > s\} E[E(H^{N,i}(t)|\mathcal{H}_{\hat{\tau}_i^N}^N)|\mathcal{H}_s^N].$$

In order to prove (A.23), we now show that

$$E(H^{N,i}(t)|\mathcal{H}_{\hat{\tau}_i^N}^N) = 0. \tag{A.25}$$

This follows informally since on the event $\{\eta_i > 0\}$, we have

$$E(H^{N,i}(t)|\mathcal{H}_{\hat{\tau}_i^N}^N) = 1\{\eta_i > 0\} \frac{E(H^{N,i}(t)|\hat{\tau}_i^N)}{P(\eta_i > 0)} = 0, \tag{A.26}$$

where the last equality is by (A.21) and the independence of $\eta_i$ and $\hat{\tau}_i^N$. In order to rigorously to prove the first equality in (A.26), we make use of Lemma 3.6 of [14].

Specifically, we show that

$$\mathcal{H}_{\hat{\tau}_i^N}^N \cap \{\eta_i > 0\} \subset (\sigma\{Q_0^N\} \vee \sigma\{\tilde{\eta}_i, i \geq 1\} \tag{A.27}$$

$$\vee\, \sigma\{\xi_r, r \geq 1, \eta_p, p \geq 1, p \neq i\} \vee \sigma\{\hat{\tau}_i^N\} \vee \mathcal{N}) \tag{A.28}$$

$$\cap \{\eta_i > 0\}.$$

Note that it is enough to check (A.27) for sets which generate $\mathcal{H}^N_{\hat{\tau}^N_i}$. However, by the definition of $\mathcal{H}^N_t$, it is not difficult to see that

$$\mathcal{H}^N_{\hat{\tau}^N_i} = \sigma\{Q^N_0\} \vee \sigma\{\tilde{\eta}_i, i \geq 1\} \vee \sigma\{\xi_r, r \geq 1\}$$

(A.29)
$$\vee \sigma\{\hat{\tau}^N_r, 1\{\eta_r = 0\}, 1\{0 < \eta_r \leq s \wedge \hat{\tau}^N_i - \hat{\tau}^N_r\},$$
$$s \geq 0, r = 1, \ldots, \hat{A}^N(\hat{\tau}^N_i)\} \vee \mathcal{N}$$

(use, for example, the argument of Brémaud [3]). Then, for $l = i, i+1, \ldots, m = 1, 2, \ldots, n = 1, 2, \ldots, s_1, s_2, \ldots, s_l > 0$ and Borel sets $A, B_1, \ldots, B_m, C_1, \ldots, C_n, D_1, \ldots, E_1, \ldots, E_l$ and $F_1, \ldots, F_l$, we have, since $\hat{A}^N(\hat{\tau}^N_i) \geq l > i$, then $\hat{\tau}^N_r = \hat{\tau}^N_i, r = i+1, \ldots, l$, that

$$\{Q^N_0 \in A\} \cap \left(\bigcap_{r=1}^{m}\{\tilde{\eta}_r \in B_r\}\right) \cap \left(\bigcap_{r=1}^{n}\{\xi_r \in C_r\}\right) \cap \{\hat{A}^N(\hat{\tau}^N_i) \geq l\}$$

$$\cap \left(\bigcap_{r=1}^{l}\{\hat{\tau}^N_r \in D_r\}\right) \cap \left(\bigcap_{r=1}^{l}\{1\{\eta_r = 0\} \in E_r\}\right)$$

$$\cap \left(\bigcap_{r=1}^{l} 1\{0 < \eta_r \leq s_r \wedge \hat{\tau}^N_i - \hat{\tau}^N_r\} \in F_r\right) \cap \{\eta_i \geq 0\}$$

(A.30) $= \{Q^N_0 \in A\} \cap \left(\bigcap_{r=1}^{m}\{\tilde{\eta}_r \in B_r\}\right) \cap \left(\bigcap_{r=1}^{n}\{\xi_r \in C_r\}\right)$

$$\cap \left(\bigcap_{r=i+1}^{l} \{\hat{\tau}^N_i = \hat{\tau}^N_r\}\right) \cap \left(\bigcap_{r=1}^{i-1}\{\hat{\tau}^N_r \in D_r\}\right) \cap \left(\bigcap_{r=i}^{l}\{\hat{\tau}^N_r \in D_r\}\right)$$

$$\cap \left(\bigcap_{r=1, r \neq i}^{l} \{1\{\eta_r = 0\} \in E_r\}\right)$$

$$\cap \left(\bigcap_{r=1}^{i-1} 1\{0 < \eta_r \leq s_r \wedge \hat{\tau}^N_i - \hat{\tau}^N_r\} \in F_r\right) \cap \{\eta_i \geq 0\},$$

when $0 \in E_i, 0 \in F_r, i \leq r \leq l$, and the left-hand side is $\varnothing$ otherwise. We show that the event on the right-hand side of (A.30) is in $(\sigma\{Q^N_0\} \vee \sigma\{\tilde{\eta}_i, i \geq 1\} \vee \sigma\{\xi_r, r \geq 1, \eta_p, p \geq 1, p \neq i\} \vee \sigma\{\hat{\tau}^N_i\} \vee \mathcal{N}) \cap \{\eta_i > 0\}$. It is enough to prove that this holds for the event

$$\bigcap_{r=i+1}^{l} \{\hat{\tau}^N_i = \hat{\tau}^N_r\} \cap \{\eta_i > 0\}.$$

First, note that for $i+1$, there exists a Borel function $h^N_{i+1}$ such that

$$\hat{\tau}^N_{i+1} = \tau_{i+1} \vee ((\hat{\tau}^N_i + \eta_i) \wedge h^N_{i+1}(Q^n_0, \tilde{\eta}_l, l \geq 1, \xi_r, r \geq 1, \eta_p, p \geq 1, p \neq i)).$$

The random variable $h_{i+1}^N(Q_0^n, \tilde{\eta}_l, l \geq 1, \xi_r, r \geq 1, \eta_p, p \geq 1, p \neq i)$ is either equal to $\tau_i$ if a server is idle in the system after the arrival of customer $i$ or, if not, the time of the next departure from the queue after the arrival of customer $i$ not including customer $i$. Since,

$$\{\hat{\tau}_{i+1}^N = \hat{\tau}_i^N\} = \{\tau_{i+1}^N \leq \hat{\tau}_i^N\} \cap \{\eta_i = 0\}$$
$$\cup \{h_{i+1}^N(Q_0^n, \tilde{\eta}_l, l \geq 1, \xi_r, r \geq 1, \eta_p, p \geq 1, p \neq i) \leq \hat{\tau}_i^N\},$$

we have that

(A.31)
$$\{\hat{\tau}_i^N = \hat{\tau}_r^N\} \cap \{\eta_i > 0\}$$
$$= \{\tau_{i+1}^N \leq \hat{\tau}_i^N\} \cap \{h_{i+1}^N(\xi_r, r \geq 1, \eta_p, p \geq 1, p \neq i) \leq \hat{\tau}_i^N\}$$
$$\cap \{\eta_i > 0\}.$$

A similar argument shows that for $r = i+2, \ldots, l$

(A.32)
$$\{\hat{\tau}_r^N = \hat{\tau}_i^N\} \cap \{\eta_i > 0\}$$
$$= \{\tau_r^N \leq \hat{\tau}_i^N\} \cap \{h_r^N(\xi_r, r \geq 1, \eta_p, p \geq 1, p \neq i) \leq \hat{\tau}_i^N\} \cap \{\eta_i > 0\},$$

where $h_r^N$ is also a Borel function and $h_r^N(\xi_r, r \geq 1, \eta_p, p \geq 1, p \neq i)$ is the time of the $(r-i)$th departure from the queue after the arrival of customer $i$ not taking into account customer $i$. By (A.31) and (A.32),

$$\bigcap_{r=i+1}^{l} \{\hat{\tau}_i^N = \hat{\tau}_r^N\} \cap \{\eta_i > 0\}$$
$$= \bigcap_{r=i+1}^{l} \{\tau_r^N \leq \hat{\tau}_i^N\} \cap \bigcap_{r=i+1}^{l} \{h_r^N(\xi_r, r \geq 1, \eta_p, p \geq 1, p \neq i) = \hat{\tau}_i^N\}$$
$$\cap \{\eta_i > 0\},$$

which yields (A.27). By Lemma 3.6 of [14], (A.27) implies the first equality in (A.26) which show that (A.25) holds and so (A.23) is proved.

We now proceed to show that (A.24) holds. We have

$$1\{\hat{\tau}_i^N \leq s\} E[H^{N,i}(t)|\mathcal{H}_s^N]$$
$$= 1\{\eta_i \leq s - \hat{\tau}_i^N\} E[H^{N,i}(t)|\mathcal{H}_s^N] + 1\{\eta_i > s - \hat{\tau}_i^N \geq 0\} E[H^{N,i}(t)|\mathcal{H}_s^N].$$

However, since both $1\{\eta_i \leq s - \hat{\tau}_i^N\}$ and $1\{\eta_i > s - \hat{\tau}_i^N\}$ are both $\mathcal{H}_s^N$-measurable, and by (A.21) we have that

$$1\{\eta_i \leq s - \hat{\tau}_i^N\} H^{N,i}(t)$$
$$= 1\{0 < \eta_i \leq s - \hat{\tau}_i^N\}$$
$$- 1\{\eta_i \leq s - \hat{\tau}_i^N\} \int_{0+}^{\eta_i \wedge (s-\hat{\tau}_i^N)^+} \frac{dF(u)}{1 - F(u-)},$$

where the latter is $\mathcal{H}_s^N$-measurable, it follows that

$$1\{\hat{\tau}_i^N \leq s\} E[H^{N,i}(t)|\mathcal{H}_s^N]$$
$$= 1\{0 < \eta_i \leq s - \hat{\tau}_i^N\}$$
(A.33)
$$- 1\{\eta_i \leq s - \hat{\tau}_i^N\} \int_{0+}^{\eta_i \wedge (s-\hat{\tau}_i^N)^+} \frac{dF(u)}{1 - F(u-)}$$
$$+ 1\{\eta_i > s - \hat{\tau}_i^N \geq 0\} E[H^{N,i}(t)|\mathcal{H}_s^N].$$

We now proceed to evaluate the quantity $1\{\eta_i > s - \hat{\tau}_i^N \geq 0\} E[H^{N,i}(t)|\mathcal{H}_s^N]$. First observe that on the event that customer $i$ has not completed service by time $s$, mathematically, the event $\{\eta_i > s - \hat{\tau}_i^N \geq 0\}$, we have that $\eta_i$ is independent of the process $\hat{A}^N$ up to time $s$ which keeps track of the number of customers who have entered service by time $s$. Also, since $\eta_i$ is independent of $\eta_l$ for $l \neq i$, we conclude from the definition of $\mathcal{H}_t^N$ and (A.21) that, on the event $\{\eta_i > s - \hat{\tau}_i^N \geq 0\}$, $H^{N,i}(t)$ is dependent upon $\mathcal{H}_s^N$ only though $\eta_i$ and $\hat{\tau}_i^N$. To put it more accurately, let $\check{A}^N(u), u \geq 0$, be the number of arrivals to the servers up until time $u$ that would have occurred if the customer with service time $\eta_i$ remained in service forever. Then, $\check{A}^N(u)$ is a Borel function of $\xi_r, r \geq 1, \eta_p, p \geq 1, p \neq i$, on the one hand, and coincides with $\hat{A}^N(u)$ for $u \leq s$ on the event $\{\hat{\tau}_i^N + \eta_i > s\}$, on the other hand. In analogy with (A.27), this yields by the definition of $\mathcal{H}_s^N$,

(A.34)
$$\mathcal{H}_s^N \cap \{\eta_i > s - \hat{\tau}_i^N \geq 0\}$$
$$\subset (\sigma\{\xi_r, r \geq 1, \eta_p, p \geq 1, p \neq i\} \vee \sigma\{\hat{\tau}_i^N\} \vee \mathcal{N})$$
$$\cap \{\eta_i > s - \hat{\tau}_i^N \geq 0\}.$$

Now noting that $1\{\eta_i > s - \hat{\tau}_i^N \geq 0\}$ is $\mathcal{H}_s^N$-measurable and applying Lemma 3.6 of [14], it then follows that

$$1\{\eta_i > s - \hat{\tau}_i^N \geq 0\} E[H^{N,i}(t)|\mathcal{H}_s^N]$$
$$= 1\{\eta_i > s - \hat{\tau}_i^N \geq 0\} \frac{E[1\{\eta_i > s - \hat{\tau}_i^N\} H^{N,i}(t)|\hat{\tau}_i^N]}{P(\eta_i > s - \hat{\tau}_i^N|\hat{\tau}_i^N)},$$

where $0/0 = 0$. Now evaluating the right-hand side of the above using (A.21), it follows

$$1\{\eta_i > s - \hat{\tau}_i^N \geq 0\} E[H^{N,i}(t)|\mathcal{H}_s^N]$$
$$= -1\{\eta_i > s - \hat{\tau}_i^N \geq 0\} \int_{0+}^{s - \hat{\tau}_i^N} \frac{dF(u)}{1 - F(u-)},$$

which, together with (A.33), implies (A.24). We now have that (A.23) and (A.24) have been proved and so the martingale property of $H^{N,i}$ has been proved.

The proof of part 2 follows identically to the proof of part 2 of Lemma 3.5 of [14]. In particular, since the second term on the right-hand side of (A.21) is $\mathbf{H}^N$-predictable, the $\mathbf{H}^N$-predictable measure of the jumps of the process $(1\{0 < \eta_i \leq t - \hat{\tau}_i^N\}, t \geq 0)$ is (see Liptser and Shiryaev [15], Jacod and Shiryaev [9])

$$\nu^{N,i}([0,T], A) = \{1 \in A\} \int_0^t 1\{\hat{\tau}_i^N < u \leq \eta_i + \hat{\tau}_i^N\} \frac{dF(u - \hat{\tau}_i^N)}{1 - F((u - \hat{\tau}_i^N)-)}$$

and so the predictable quadratic-variation process of $H^{N,i}$ is (see, e.g., Liptser and Shiryaev [15], Problem 11, Chapter 4, Section 1)

$$\langle H^{N,i}\rangle(t) = \int_0^t \int_{\mathbb{R}} x^2 \nu^{N,i}(du, dx) - \sum_{0 < u \leq t} \left( \int_{\mathbb{R}} x \nu^{N,i}(\{u\}, dx) \right)^2$$

$$= \int_0^t 1\{\hat{\tau}_i^N < u \leq \eta_i + \hat{\tau}_i^N\} \frac{dF(u - \hat{\tau}_i^N)}{1 - F((u - \hat{\tau}_i^N)-)}$$

$$- \sum_{0 < u \leq t} 1\{\hat{\tau}_i^N < u \leq \eta_i + \hat{\tau}_i^N\} \left( \frac{\Delta F(u - \hat{\tau}_i^N)}{1 - F((u - \hat{\tau}_i^N)-)} \right)^2$$

$$= \int_0^{\eta_i \vee (t - \hat{\tau}_i^N)^+} \frac{dF(u)}{1 - F(u-)} - \sum_{0 < u \leq \eta_i \wedge (t - \hat{\tau}_i^N)^+} \left( \frac{\Delta F(u)}{1 - F(u-)} \right)^2,$$

where the sum is over all of the jumps. Since the above is equal to the right-hand side of (A.22), this then completes the proof of part 2.

We now demonstrate step 3 by proving the martingale property for $H^{N,i}H^{N,j}$ in a similar to manner to the proof of the martingale property for $H^{N,i}$. Specifically, for $s < t, j < i$, we prove that

(A.35) $$1\{\hat{\tau}_i^N > s\} E[H^{N,i}(t) H^{N,j}(t) | \mathcal{H}_s^N] = 0$$

and

(A.36) $$1\{\hat{\tau}_i^N \leq s\} E[H^{N,i}(t) H^{N,j}(t) | \mathcal{H}_s^N] = H^{N,i}(t) H^{N,j}(t).$$

For (A.35), we have

(A.37) $$1\{\hat{\tau}_i^N > s\} E[H^{N,i}(t) H^{N,j}(t) | \mathcal{H}_s^N]$$

(A.38) $$= 1\{\hat{\tau}_i^N > s\} E[E[H^{N,i}(t) H^{N,j}(t) | H_{\hat{\tau}_i^N}^N] | \mathcal{H}_s^N].$$

However, since

(A.39) $$E[H^{N,i}(t) H^{N,j}(t) | H_{\hat{\tau}_i^N}^N]$$
$$= 1\{\eta_j \leq \hat{\tau}_i^N - \hat{\tau}_j^N\} E[H^{N,i}(t) H^{N,j}(t) | H_{\hat{\tau}_i^N}^N]$$
$$+ 1\{\eta_j > \hat{\tau}_i^N - \hat{\tau}_j^N\} E[H^{N,i}(t) H^{N,j}(t) | H_{\hat{\tau}_i^N}^N],$$

and $1\{\eta_j \leq \hat{\tau}_i^N - \hat{\tau}_j^N\}$ and $1\{\eta_j \leq \hat{\tau}_i^N - \hat{\tau}_j^N\}H^{N,j}(t) = 1\{\eta_j \leq \hat{\tau}_i^N - \hat{\tau}_j^N\}H^{N,j}(t \wedge \hat{\tau}_i^N)$ are $\mathcal{H}_{\hat{\tau}_i^N}^N$-measurable [use (A.29)], it follows that

$$1\{\eta_j \leq \hat{\tau}_i^N - \hat{\tau}_j^N\}E[H^{N,i}(t)H^{N,j}(t)|H_{\hat{\tau}_i^N}^N]$$
$$= E[1\{\eta_j \leq \hat{\tau}_i^N - \hat{\tau}_j^N\}H^{N,i}(t)H^{N,j}(t)|H_{\hat{\tau}_i^N}^N]$$
$$= 1\{\eta_j \leq \hat{\tau}_i^N - \hat{\tau}_j^N\}H^{N,j}(t)E[H^{N,i}(t)|H_{\hat{\tau}_i^N}^N],$$

and so, since $H^{N,i}$ is a square-integrable martingale, by Doob's stopping theorem ([9], I.1.19, I.1.42), $E[H^{N,i}(t)|H_{\hat{\tau}_i^N}^N] = H^{N,i}(\hat{\tau}_i^N) = 0$ and thus, the first term on the right-hand side of (A.39) is 0.

Consider now the second term. On the event $\{\eta_j > \hat{\tau}_i^N - \hat{\tau}_j^N\}$, we have that customer $j$ finishes service after customer $i$ arrives and so customer $j$'s service time has no no effect on $\hat{\tau}_i^N$, the time at which customer $i$ enters service. Thus, $\eta_j$ and $\hat{\tau}_i^N$ are independent on the event $\{\eta_j > \hat{\tau}_i^N - \hat{\tau}_j^N\}$. To put it more accurately, there exists a random variable $\check{\tau}_i^N$ which is a Borel function of $\xi_r, r \geq 1, \eta_p, p \geq 1, p \neq i, p \neq j$, such that $\{\eta_j > \hat{\tau}_i^N - \hat{\tau}_j^N\} = \{\eta_j > \check{\tau}_i^N - \hat{\tau}_j^N\}$ and $\hat{\tau}_i^N = \check{\tau}_i^N$ on either event. One may view $\check{\tau}_i^N$ as the time at which customer $i$ would enters service if customer $j$'s service time were infinitely long. Thus, applying Lemma 3.6 of [14] and using the fact that $\eta_i$ and $\eta_j$ are independent of $\check{\tau}_i^N$ and $\hat{\tau}_j^N$, we have that

$$1\{\eta_j > \hat{\tau}_i^N - \hat{\tau}_j^N\}E[H^{N,i}(t)H^{N,j}(t)|H_{\hat{\tau}_i^N}^N]$$
$$= 1\{\eta_j > \check{\tau}_i^N - \hat{\tau}_j^N\}\frac{E[1\{\eta_j > \check{\tau}_i^N - \hat{\tau}_j^N\}\check{H}^{N,i}(t)H^{N,j}(t)|\check{\tau}_i^N, \hat{\tau}_j^N]}{P(\eta_j > \check{\tau}_i^N - \hat{\tau}_j^N|\check{\tau}_i^N, \hat{\tau}_j^N)},$$

where $\check{H}^{N,i}$ denotes $H^{N,i}$ with $\check{\tau}_i^N$ substituted for $\hat{\tau}_i^N$. Furthermore, since $\eta_i$ is independent of $\check{\tau}_i^N, \eta_j$ and $\hat{\tau}_j^N$, we obtain that

$$E[1\{\eta_j > \check{\tau}_i^N - \hat{\tau}_j^N\}\check{H}^{N,i}(t)H^{N,j}(t)|\check{\tau}_i^N, \hat{\tau}_j^N]$$
$$= E[1\{\eta_j > \check{\tau}_i^N - \hat{\tau}_j^N\}H^{N,j}(t)|\check{\tau}_i^N, \hat{\tau}_j^N]E[\check{H}^{N,i}(t)|\check{\tau}_i^N],$$

where the last multiplier on the right-hand side is equal to 0 by the definition of $H^{N,i}$ and the fact that $\eta_i$ is independent of $\check{\tau}_i^N$. Thus, the right-hand side of (A.39) is 0, and so $E[H^{N,i}(t)H^{N,j}(t)|H_{\hat{\tau}_i^N}^N] = 0$ so that (A.35) is proved.

In order to prove (A.36), we proceed similarly to (A.24) and so some of the details are omitted. First, note that

$$1\{\hat{\tau}_i^N \leq s\}E[H^{N,i}(t)H^{N,j}(t)|\mathcal{H}_s^N]$$
$$= 1\{\eta_i \leq s - \hat{\tau}_i^N\}1\{\eta_j \leq s - \hat{\tau}_j^N\}E[H^{N,i}(t)H^{N,j}(t)|\mathcal{H}_s^N]$$

$$\begin{aligned}
&+ 1\{\eta_i > s - \hat{\tau}_i^N \geq 0\} 1\{\eta_j \leq s - \hat{\tau}_j^N\} E[H^{N,i}(t) H^{N,j}(t) | \mathcal{H}_s^N] \\
&+ 1\{\eta_i > s - \hat{\tau}_i^N \geq 0\} 1\{\eta_j \leq s - \hat{\tau}_j^N\} E[H^{N,i}(t) H^{N,j}(t) | \mathcal{H}_s^N] \\
&+ 1\{\eta_i > s - \hat{\tau}_i^N \geq 0\} 1\{\eta_j > s - \hat{\tau}_j^N \geq 0\} \\
&\quad \times E[H^{N,i}(t) H^{N,j}(t) | \mathcal{H}_s^N]
\end{aligned} \tag{A.40}$$

and, further,

$$\begin{aligned}
&1\{\eta_i \leq s - \hat{\tau}_i^N\} 1\{\eta_j \leq s - \hat{\tau}_j^N\} E[H^{N,i}(t) H^{N,j}(t) | \mathcal{H}_s^N] \\
&= 1\{\eta_i \leq s - \hat{\tau}_i^N\} 1\{\eta_j \leq s - \hat{\tau}_j^N\} \\
&\quad \times E[1\{\eta_i \leq s - \hat{\tau}_i^N\} H^{N,i}(t) 1\{\eta_j \leq s - \hat{\tau}_j^N\} H^{N,j}(t) | \mathcal{H}_s^N] \\
&= 1\{\eta_i \leq s - \hat{\tau}_i^N\} 1\{\eta_j \leq s - \hat{\tau}_j^N\} H^{N,i}(t) H^{N,j}(t),
\end{aligned} \tag{A.41}$$

$$\begin{aligned}
&1\{\eta_i \leq s - \hat{\tau}_i^N\} 1\{\eta_j > s - \hat{\tau}_j^N \geq 0\} E[H^{N,i}(t) H^{N,j}(t) | \mathcal{H}_s^N] \\
&= 1\{\eta_i \leq s - \hat{\tau}_i^N\} H^{N,i}(s) 1\{\eta_j > s - \hat{\tau}_j^N \geq 0\} E[H^{N,j}(t) | \mathcal{H}_s^N] \\
&= 1\{\eta_i \leq s - \hat{\tau}_i^N\} H^{N,i}(s) 1\{\eta_j > s - \hat{\tau}_j^N \geq 0\} H^{N,j}(s),
\end{aligned} \tag{A.42}$$

$$\begin{aligned}
&1\{\eta_i > s - \hat{\tau}_i^N \geq 0\} 1\{\eta_j \leq s - \hat{\tau}_j^N\} E[H^{N,i}(t) H^{N,j}(t) | \mathcal{H}_s^N] \\
&= 1\{\eta_j \leq s - \hat{\tau}_j^N 0\} H^{N,j}(s) 1\{\eta_i > s - \hat{\tau}_i^N \geq 0\} E[H^{N,i}(t) | \mathcal{H}_s^N] \\
&= 1\{\eta_j \leq s - \hat{\tau}_j^N 0\} H^{N,j}(s) 1\{\eta_i > s - \hat{\tau}_i^N \geq 0\} H^{N,i}(s),
\end{aligned} \tag{A.43}$$

where in (A.42) and (A.43) we use the martingale property of $H^{N,i}$ and $H^{N,j}$.

Consider now the last term on the right of (A.40). Since $\{s - \hat{\tau}_i^N \geq 0, \eta_j > s - \hat{\tau}_j^N\} \subset \{\eta_j > \hat{\tau}_i^N - \hat{\tau}_j^N\}$, it follows as above that $\{s - \hat{\tau}_i^N \geq 0, \eta_j > s - \hat{\tau}_j^N\} = \{s - \check{\tau}_i^N \geq 0, \eta_j > s - \hat{\tau}_j^N\}$ and $\hat{\tau}_i^N = \check{\tau}_i^N$ on either event. Hence,

$$\begin{aligned}
&1\{\eta_i > s - \hat{\tau}_i^N \geq 0\} 1\{\eta_j > s - \hat{\tau}_j^N \geq 0\} E[H^{N,i}(t) H^{N,j}(t) | \mathcal{H}_s^N] \\
&= 1\{\eta_i > s - \check{\tau}_i^N \geq 0\} 1\{\eta_j > s - \hat{\tau}_j^N \geq 0\} E[\check{H}^{N,i}(t) H^{N,j}(t) | \mathcal{H}_s^N],
\end{aligned}$$

where $\eta_i$ and $\eta_j$ are independent of $\check{\tau}_i^N$ and $\hat{\tau}_j^N$. Also, similar to (A.27) and (A.34),

$$\begin{aligned}
&\mathcal{H}_s^N \cap \{\eta_i > s - \check{\tau}_i^N\} \cap \{\eta_j > s - \hat{\tau}_j^N\} \\
&\quad \subset (\sigma\{\xi_r, r \geq 1, \eta_p, p \geq 1, p \neq i, p \neq j\} \vee \sigma\{\check{\tau}_i^N, \hat{\tau}_j^N\} \vee \mathcal{N}) \\
&\qquad \cap \{\eta_i > s - \check{\tau}_i^N\} \cap \{\eta_j > s - \hat{\tau}_j^N\},
\end{aligned}$$

and so by Lemma 3.6 of [14],

$$1\{\eta_i > s - \hat{\tau}_i^N \geq 0\}1\{\eta_j > s - \hat{\tau}_j^N \geq 0\}E[H^{N,i}(t)H^{N,j}(t)|\mathcal{H}_s^N]$$
$$= 1\{\eta_i > s - \hat{\tau}_i^N \geq 0\}1\{\eta_j > s - \hat{\tau}_j^N \geq 0\}$$
(A.44)
$$\times (E[1\{\eta_i > s - \check{\tau}_i^N \geq 0\}1\{\eta_j > s - \hat{\tau}_j^N \geq 0\}\check{H}^{N,i}(t)$$
$$\times H^{N,j}(t)|\check{\tau}_i^N, \hat{\tau}_j^N])$$
$$/(P(\eta_i > s - \check{\tau}_i^N \geq 0, \eta_j > s - \hat{\tau}_j^N \geq 0|\check{\tau}_i^N, \hat{\tau}_j^N)).$$

However, since $\eta_j$ and $\eta_i$ are independent of each other and $\check{\tau}_i^N$ and $\hat{\tau}_j^N$, we have that

$$E[1\{\eta_i > s - \check{\tau}_i^N \geq 0\}1\{\eta_j > s - \hat{\tau}_j^N \geq 0\}\check{H}^{N,i}(t)H^{N,j}(t)|\check{\tau}_i^N, \hat{\tau}_j^N]$$
(A.45)
$$= E[1\{\eta_i > s - \check{\tau}_i^N \geq 0\}\check{H}^{N,i}(t)|\check{\tau}_i^N]$$
$$\times E[1\{\eta_j > s - \hat{\tau}_j^N \geq 0\}H^{N,j}(t)|\hat{\tau}_j^N]$$

and

(A.46)
$$P(\eta_i > s - \check{\tau}_i^N \geq 0, \eta_j > s - \hat{\tau}_j^N \geq 0|\check{\tau}_i^N, \hat{\tau}_j^N)$$
$$= P(\eta_i > s - \check{\tau}_i^N \geq 0|\check{\tau}_i^N)P(\eta_j > s - \hat{\tau}_j^N \geq 0|\hat{\tau}_j^N).$$

Applying Lemma 3.6 of [14] and using analogues of (A.34), we have that

$$1\{\eta_i > s - \check{\tau}_i^N \geq 0\}\frac{E[1\{\eta_i > s - \check{\tau}_i^N \geq 0\}\check{H}^{N,i}(t)|\check{\tau}_i^N]}{P(\eta_i > s - \check{\tau}_i^N \geq 0|\check{\tau}_i^N)}$$
$$= 1\{\eta_i > s - \check{\tau}_i^N \geq 0\}E[\check{H}^{N,i}(t)|\mathcal{H}_s^N]$$

and

$$1\{\eta_j > s - \hat{\tau}_j^N \geq 0\}\frac{E[1\{\eta_j > s - \hat{\tau}_j^N \geq 0\}H^{N,j}(t)|\hat{\tau}_j^N]}{P(\eta_j > s - \hat{\tau}_j^N \geq 0|\hat{\tau}_j^N)}$$
$$= 1\{\eta_j > s - \hat{\tau}_j^N \geq 0\}E[H^{N,j}(t)|\mathcal{H}_s^N],$$

so that dividing (A.45) by (A.46), using (A.44), the fact that $\hat{\tau}_i^N = \check{\tau}_i^N$ on $\{s - \hat{\tau}_i^N \geq 0, \eta_j > s - \hat{\tau}_j^N\} = \{s - \check{\tau}_i^N \geq 0, \eta_j > s - \hat{\tau}_j^N\}$ and the martingale property of $H^{N,i}$ and $H^{N,j}$ we get

$$1\{\eta_i > s - \hat{\tau}_i^N \geq 0\}1\{\eta_j > s - \hat{\tau}_j^N \geq 0\}E[H^{N,i}(t)H^{N,j}(t)|\mathcal{H}_s^N]$$
$$= 1\{\eta_i > s - \hat{\tau}_i^N \geq 0\}E[H^{N,i}(t)|\mathcal{H}_s^N]$$
(A.47)
$$\times 1\{\eta_j > s - \hat{\tau}_j^N \geq 0\}E[H^{N,j}(t)|\mathcal{H}_s^N]$$
$$= 1\{\eta_i > s - \hat{\tau}_i^N \geq 0\}H^{N,i}(s)1\{\eta_j > s - \hat{\tau}_j^N \geq 0\}H^{N,j}(s).$$

Substituting (A.41)–(A.43) and (A.47) into (A.40), we obtain (A.36), which completes the proof of the lemma. □

Now note that by (4.3) and (A.19), we have

$$\bar{M}_2^N(t) = \bar{G}^N(t) + \bar{H}^N(t), \tag{A.48}$$

where

$$\bar{G}^N = \frac{G^N}{N} \tag{A.49}$$

and

$$\bar{H}^N = \frac{H^N}{N}. \tag{A.50}$$

It therefore follows by (A.48) that in order to show $\bar{M}_2^N \Rightarrow 0$ as $N \to \infty$, it will be sufficient to show that $\bar{G}^N \Rightarrow 0$ and $\bar{H}^N \Rightarrow 0$ as $N \to \infty$. First, however, we must provide the following result.

Let

$$\check{A}^N(t) = \frac{\hat{A}^N(t)}{N}, \tag{A.51}$$

be the fluid scaled number of customers to enter service by time $t$. We then have the following.

LEMMA A.2. *For each $T \geq 0$, there exists a $\kappa \geq 0$ such that $P(\check{A}^N(T) \geq \kappa) \to 0$ as $N \to \infty$.*

PROOF. In order to show that the result is true, we stochastically bound $\{\check{A}^N(T)\}$ by another sequence of random variables for which the result holds. This will then imply that the result holds for $\{\check{A}^N(T)\}$ as well.

Let $\min(Q^N(T), N)$ be the total number of customers in service at time $T$. We then have that

$$\min(Q^N(T), N) = \min(Q_0^N, N) + \hat{A}^N(T) - D^N(T), \tag{A.52}$$

where $D^N(T)$ is the number of departures from the system by time $T$. Equation (A.52) then implies that

$$\begin{aligned}
\hat{A}^N(T) &= \min(Q^N(T), N) + D^N(T) - \min(Q_0^N, N) \\
&\leq \min(Q^N(T), N) + D^N(T) \\
&\leq N + D^N(T).
\end{aligned} \tag{A.53}$$

We next bound $D^N(T)$. Let $S_i^N(t)$ be the number of departures from server $i$ in its first $t$ units of processing time for $t \geq 0$ and let $B_i^N(t)$ be the amount of time that server $i$ is busy in the first $t$ time units. We then have that

$$D^N(T) = \sum_{i=1}^N S_i^N(B_i^N(T)) \leq \sum_{i=1}^N S_i^N(T),$$

since $B_i^N(T) \leq T$.

Now note that for each $i$, $S_i^N(T)$ is either the number of renewals by time $T$ of a pure renewal process with interarrival distribution $F$ or a delayed renewal process with delay distribution $F_0$ and interarrival distribution $F$. Furthermore, for $i \neq j$, we have that $S_i^N(T)$ and $S_j^N(T)$ are independent of one another. Letting $\{P_i, i \geq 1\}$ be an i.i.d. sequence of pure renewal processes with interarrival distribution $F$ and $\{Q_i, i \geq 1\}$ an i.i.d. sequence of delayed renewal processes with delay distribution $F_0$ and interarrival distribution $F$, it therefore follows that

$$N^{-1} \sum_{i=1}^N S_i^N(T) \leq^{st} N^{-1} \sum_{i=1}^N P_i(T) + N^{-1} \sum_{i=1}^N Q_i(T)$$
$$\Rightarrow M(T) + M_D(T) \quad \text{as } N \to \infty,$$

where $M$ is the renewal function associated with $P_1$ and $M_D$ is the renewal function associated with $Q_1$. This completes the proof. □

We now show that $\bar{M}_2^N \Rightarrow 0$ as $N \to \infty$. We begin by showing that $\bar{G}^N \Rightarrow 0$ as $N \to \infty$. Let

(A.54) $$\bar{U}^N = \frac{U^N}{N}.$$

We then have the following.

LEMMA A.3. $\bar{G}^N \Rightarrow 0$ as $N \to \infty$.

PROOF. The proof of is nearly identical to the proof of Lemma 3.4 of [14] but for completeness we will include it here as well.

We first show that for each $\delta > 0$ and $T > 0$,

$$\lim_{\varepsilon \downarrow 0} \limsup_N P\left(\sup_{0 \leq t \leq T} \left| \int_0^t \frac{\bar{V}^N(t-x, x-)}{1 - F(x-)} \right.\right.$$

(A.55) $$\left.\left. \times 1\{F(x-) > 1 - \varepsilon\} \, dF(x) \right| > \delta\right)$$
$$= 0,$$

where $\bar{V}^N(t, x) = N^{-1} V^N(t, x)$.

By (A.14) and recalling the definition of $\bar{U}^N$ from (A.54), we have that for any $k > 0$,

$$P\left(\sup_{0 \leq t \leq T} \left| \int_0^t \frac{\bar{V}^N(t-x, x-)}{1 - F(x-)} 1\{F(x-) > 1 - \varepsilon\} \, dF(x) \right| > \delta \right)$$
$$\leq P(\check{A}^N(T) > kT)$$
$$+ P\left( \int_0^\infty \frac{1\{F(x-) > 1 - \varepsilon\}}{1 - F(x-)} \sup_{0 \leq t \leq kT} |\bar{U}^N(t, F(x-))| \, dF(x) > \delta \right).$$

For $k$ sufficiently large, we have by Lemma A.2 that

$$P(\check{A}^N(T) > kT) \to 0 \quad \text{as } N \to \infty.$$

Thus, by applying Chebyshev's inequality and Fubini's theorem, we now must prove that

$$\lim_{\varepsilon \downarrow 0} \limsup_N \int_0^\infty \frac{1\{F(x-) > 1 - \varepsilon\}}{1 - F(x-)} E \sup_{0 \leq t \leq kT} |\bar{U}^N(t, F(x-))| \, dF(x) = 0.$$

However, the proof of this proceeds identically to as in Lemma 3.4 of [14], which completes the proof. □

We next show that $\bar{H}^N$ converges to 0 as $N$ goes to $\infty$. Again, the modifications to the proof of Lemma 3.7 of [14] are slight but we include a full proof for completeness.

LEMMA A.4. $\bar{H}^N \Rightarrow 0$ as $N \to \infty$.

PROOF. Let

$$\hat{H}^N(t) = N^{-1} \sum_{i=1}^{\hat{A}^N(t)} \left( 1\{0 < \eta_i \leq t - \hat{\tau}_i^N\} - \int_{0+}^{\eta_i \wedge (t - \hat{\tau}_i^N)^+} \frac{dF(u)}{1 - F(u-)} \right), \quad t \geq 0,$$

and note that by (A.17), (A.20) and (A.50) we have that

$$\bar{H}^N(t) = N^{-1} \sum_{i=1}^{\hat{A}^N(t)} (1\{\eta_i = 0\} - F(0)) + \hat{H}^N(t).$$

We first show that the term involving the summation converges to 0. Let $T \geq 0$ and $\delta > 0$. We have

$$P\left( \sup_{0 \leq t \leq T} \left| N^{-1} \sum_{i=1}^{\hat{A}^N(t)} (1\{\eta_i = 0\} - F(0)) \right| > \delta \right)$$
$$\leq P(N^{-1} \hat{A}^N(T) > k) + P\left( \sup_{0 \leq t \leq 1} \left| N^{-1} \sum_{i=1}^{\lfloor Nkt \rfloor} (1\{\eta_i = 0\} - F(0)) \right| > \delta \right).$$

However, for sufficiently large $k$, we have by Lemma A.2 that $P(N^{-1}\hat{A}^N(T) > k) \to 0$ as $N \to \infty$. Furthermore, by the functional strong law of large numbers and the i.i.d. assumption of $\{\eta_i, i \geq 1\}$, it follows that

$$P\left(\sup_{0 \leq t \leq 1} \left| N^{-1} \sum_{i=1}^{\lfloor Nkt \rfloor} (1\{\eta_i = 0\} - F(0)) \right| > \delta \right) \to 0 \qquad \text{as } N \to \infty.$$

It thus remains to show the convergence of $\hat{H}^N$ to 0.

Fix $T > 0$. For each $\varepsilon > 0$, we have

$$P\left(\sup_{0 \leq t \leq T} \hat{H}^N(t) > \varepsilon\right) \leq P(N^{-1}\hat{A}^N(T) > k) + P\left(\sup_{0 \leq t \leq T} |\bar{H}^N_{Nk}(t)| > \varepsilon\right),$$

where

$$\bar{H}^N_{Nk} = \frac{H^N_{Nk}}{N}.$$

By definition (A.15) and Lemma A.2, for $k$ sufficiently large,

$$P(N^{-1}\hat{A}^N(T) > k) \to 0 \qquad \text{as } N \to \infty.$$

Next, recall by Lemma A.1 that $\bar{H}^N_{Nk}$ is an $\mathbf{H}^N$-square-integrable martingale with predictable quadratic variation process

(A.56)
$$\langle \bar{H}^N_{Nk} \rangle(t)$$
$$= N^{-2} \sum_{i=1}^{\hat{A}^N(t) \wedge Nk} \int_{0+}^{\eta_i \wedge (t-\hat{\tau}^N_i)^+} \frac{1 - F(u)}{(1 - F(u-))^2} dF(u), \qquad t \geq 0.$$

Thus, by the Lenglart–Rebolledo inequality [15], it follows that for any $\gamma > 0$,

$$P\left(\sup_{0 \leq t \leq T} |\bar{H}^N_{Nk}(t)| > \varepsilon\right) \leq \frac{\gamma}{\varepsilon^2} + P(\langle \bar{H}^N_{Nk} \rangle(T) > \gamma).$$

However, by (A.56),

(A.57) $$\langle \bar{H}^N_{Nk} \rangle(T) \leq N^{-2} \sum_{i=1}^{\hat{A}^N(T)} \int_0^{\eta_i} \frac{dF(u)}{1 - F(u-)}.$$

Furthermore, since $E[\int_0^{\eta_i} (1 - F(u-))^{-1} dF(u)] = 1$, it follows by the functional strong law of large numbers that

$$N^{-2} \sum_{i=1}^{\lfloor N \cdot \rfloor} \int_0^{\eta_i} \frac{dF(u)}{1 - F(u-)} \Rightarrow 0 \qquad \text{as } N \to \infty.$$

By (A.57), the Random Time Change theorem [1] and Lemma A.2, this then implies that for any $\gamma > 0$,

$$P(\langle \bar{H}^N_{Nk} \rangle(T) > \gamma) \to 0 \qquad \text{as } N \to \infty,$$

which completes the proof. □

We are now in a position to give a proof of Proposition 4.2.

PROOF OF PROPOSITION 4.2. The proof follows by the decomposition (A.48) and Lemmas A.3 and A.4 above. □

We now proceed to proving Proposition 5.1. We first begin with the following result. Recall the definition of $\hat{A}^N(t)$ as the number of customers in the $N$th system who entered service after time $0-$ but before or at time $t$. Also, recall the definition of $\check{A}^N$ from (A.51) as the fluid scaled version of $\hat{A}^N$. We then have the following result.

LEMMA A.5. *Under the assumptions of Section 5, $\check{A}^N \Rightarrow e$ as $N \to \infty$.*

PROOF. First, note the relationship
$$\hat{A}^N(t) = A^N(t) - (Q^N(t) - N)^+ + (Q^N(0) - N)^+, \qquad t \geq 0,$$
which, dividing by $N$, may be equivalently expressed as
$$\check{A}^N(t) = \bar{A}^N(t) - (\bar{Q}^N(t) - 1)^+ + (\bar{Q}_0^N - 1)^+, \qquad t \geq 0.$$
By (5.1), (5.2) and the Halfin–Whitt assumption (5.3), it follows that $\bar{A}^N \Rightarrow e$ as $N \to \infty$. Thus, by Corollary 5.1, the assumption that $\bar{Q}_0^N \Rightarrow 1$ as $N \to \infty$ and the Continuous Mapping theorem [20], it follows that
$$\check{A}^N = \bar{A}^N - (\bar{Q}^N - 1)^+ + (\bar{Q}_0^N - 1)^+ \Rightarrow e \qquad \text{as } N \to \infty,$$
which completes the proof. □

Now, define the processes
$$\tilde{G}^N = \frac{G^N}{\sqrt{N}}$$
and
$$\tilde{H}^N = \frac{H^N}{\sqrt{N}},$$
and note that by (5.9) and (A.19) it follows that
$$\tilde{M}_2^N = \tilde{G}^N + \tilde{H}^N. \tag{A.58}$$

Our next result will be to show that the sequence $\{\tilde{M}_2^N\}$ is tight. In order to show this, it will be sufficient to show that both $\{\tilde{G}^N\}$ and $\{\tilde{H}^N\}$ are tight. The proofs of these results are similar to proofs of Lemmas 3.4 and 3.7 of [14], and hence have not been included. We begin with a proof of the tightness of $\{\tilde{G}^N\}$.

LEMMA A.6. *The sequence $\{\tilde{G}^N\}$ is tight.*

PROOF. By virtue of Lemma A.5 and the fact that the identity process $e(t) = t$ is a continuous process, the proof now follows identically to the proof of Lemma 3.4 in [14]. The modifications to this proof are essentially trivial and the interested reader is referred to Lemma 3.4 of [14] for further details. □

Next, we show that $\{\tilde{H}^N\}$ is tight.

LEMMA A.7. *The sequence $\{\tilde{H}^N\}$ is tight.*

PROOF. Since by Lemma A.1, the process $H_k^N$ is an $\mathbf{H}^N$-square-integrable-martingale for each $N$ and $k$, by Lemma A.5 the proof now follows similarly to the proof of Lemma 3.7 of [14] and will not be included. Again, the interested reader is referred to [14] for further details. □

We may now state the following result.

PROPOSITION A.1. *The sequence $\{\tilde{M}_2^N\}$ is tight.*

PROOF. The result follows by the decomposition (A.58) and Lemmas A.6 and A.7 above. □

We are now ready to give a proof of Proposition 5.1. Before doing so, however, we must first recall Lemma 5.2 from [14]. The proof of this result is similar in our case and therefore will not be included for the sake of brevity.

Let $\beta_i(x, y)$ be bounded real-valued Borel functions such that $E[\beta_i(x, \eta_i) = 0]$ and define the processes by $R_m^N = \{R_m^N(t), t \geq 0\}$ and $\langle R_m^N \rangle = \{\langle R_m^N \rangle(t), t \geq 0\}$, $m = 1, 2, \ldots$, by

$$(A.59) \quad R_m^N(t) = \sum_{i=1}^{\hat{A}^N(t) \wedge m} \beta_i(\hat{\tau}_i^N, \eta_i) \quad \text{and} \quad \langle R_m^N \rangle(t) = \sum_{i=1}^{\hat{A}^N(t) \wedge m} \bar{\beta}_i(\hat{\tau}_i^N),$$

where

$$\bar{\beta}_i(x) = E\beta_i^2(x, \eta_i).$$

We also set the $\sigma$-fields $\hat{\mathcal{F}}_t^N = \sigma\{\hat{\tau}_i^N, \eta_i, 1 \leq i \leq \lfloor t \rfloor\} \vee \mathcal{N}$ and $\mathcal{F}_t^N = \sigma\{\hat{\tau}_i^N \wedge \hat{\tau}_{\hat{A}^N(t)+1}^N, \eta_{i \wedge \hat{A}^N(t)}, i \geq 1\} \vee \mathcal{N}$, and define the filtrations $\hat{\mathbf{F}}^N = \{\hat{\mathcal{F}}_t^N, t \geq 0\}$ and $\mathbf{F}^N = \{\mathcal{F}_t^N, t \geq 0\}$.

We then have the following result.

LEMMA A.8. *1. The $\hat{\tau}_i^N, i=1,2,\ldots,$ are $\mathbf{F}^N$-stopping times, and the following inclusions hold: $\mathcal{F}_{\hat{\tau}_i^N}^N \supset \hat{\mathcal{F}}_{i+1}^N, \mathcal{G}_i^N \subset \hat{\mathcal{F}}_i^N$, where $\mathcal{G}_i^N = \sigma\{\mathcal{B} \cap \{\hat{\tau}_i^N > t\}, t \geq 0, \mathcal{B} \in \mathcal{F}_t^N\}$;*

*2. The process $\hat{A}^N$ is $\mathcal{F}^N$-predictable;*

*3. The processes $R_m^N, m = 1, 2, \ldots,$ are $\mathbf{F}^N$-square-integrable martingales with the processes $\langle R_m^N \rangle$ as predictable quadratic-variation processes.*

PROOF. The proof is identical to the proof of Lemma 5.2 of [14]. □

We are now prepared to give a proof of Proposition 5.1.

PROOF OD PROPOSITION 5.1. Ourproof is similar to the proof of Lemma 5.3 of [14] but we restate it here for the sake of completeness. Our first step is to show that the finite-dimensional distributions of $(\tilde{M}_2^N, \hat{M}_2^N)$ converge to those of $(\tilde{M}_2, \tilde{M}_2)$. We denote finite-dimensional convergence by $\stackrel{\text{f.d.}}{\Rightarrow}$.

Let
$$\tilde{U}^N = \frac{U^N}{\sqrt{N}}$$
and note that by Lemma 3.1 of [14], $\tilde{U}^N \Rightarrow \tilde{U}$ in $D([0,\infty), D[0,1])$ as $N \to \infty$, where $\tilde{U}$ is the Kiefer process. Next, let

$$\tilde{M}_{2,k}^N(t) = \sum_{i=1}^k \Box \tilde{U}^N((\hat{A}^N(s_{i-1}^k), F(0)), (\hat{A}^N(s_i^k), F(t-s_i^k))), \tag{A.60}$$

where the increment
$$\Box \tilde{U}^N((a_1, a_2), (b_1, b_2)) = \tilde{U}^N(b_1, b_2) - \tilde{U}^N(a_1, b_2) \\ - \tilde{U}^N(b_1, a_2) + \tilde{U}^N(a_1, a_2),$$

and the points $0 = s_0^k < s_1^k < \cdots < s_k^k = t$ are chosen such that
$$\max_{1 \leq i \leq k} |s_i^k - s_{i-1}^k| \to 0 \quad \text{as } k \to \infty.$$

We also define in analogy,
$$M_{2,k}(t) = \sum_{i=1}^K (\Box \tilde{U}((e(s_{i-1}^k), F(0)), (e(s_i^k), F(t-s_i^k)))) \\ + (\tilde{U}(e(s_i^k), F(0)) - \tilde{U}(e(s_{i-1}^k), F(0)))),$$

where
$$\Box \tilde{U}((a_1, a_2), (b_1, b_2)) = \tilde{U}(b_1, b_2) - \tilde{U}(a_1, b_2) - \tilde{U}(b_1, a_2) + \tilde{U}(a_1, a_2).$$

We now show that

(a) $\tilde{M}_{2,k}^N \stackrel{\text{f.d.}}{\Rightarrow} M_{2,k}$,
(b) $\lim_{k\to\infty} \limsup_{N\to\infty} P(|\tilde{M}_{2,k}^N(t) - \tilde{M}_2^N(t)| > \eta) = 0$ for $\eta > 0, t > 0$,
(c) $\lim_{n\to\infty} P(|\hat{M}_2^N(t) - \tilde{M}_2^N(t)| > \eta) = 0$ for $\eta > 0, t > 0$.

Since $M_{2,k}(t) \stackrel{P}{\Rightarrow} M_2(t)$ as $k \to \infty$ by definition, this will prove the finite-dimensional convergence stated in the paragraph above.

The proofs of (a) and (b) are identical to the proofs in Lemma 5.3 of [14] but we include them here for the sake of completeness. We proceed as follows.

By the Lemma 3.1 of [14] and the continuity of the Keifer process $\tilde{U}$, it follows, setting

$$(\text{A.61}) \quad \check{M}_{2,k}^N(t) = \sum_{i=1}^{k} \Box \tilde{U}^N((e(s_i^K), F(0)), (e(s_i^k), F(t - s_i^k))), \qquad t \geq 0,$$

that

$$\check{M}_{2,k}^N \Rightarrow M_{2,k} \qquad \text{as } N \to \infty.$$

Next, by Lemma A.5, Lemma 3.1 of [14] and the continuity of $\tilde{U}$ and $e$, we obtain from (A.60) and (A.61) that

$$(\text{A.62}) \quad \lim_{N\to\infty} P\Big(\sup_{0\leq t\leq T} |\tilde{M}_{2,k}^N(t) - \check{M}_{2,k}^N(t)| > \varepsilon\Big) = 0, \qquad T > 0, \varepsilon > 0.$$

This then implies that $\tilde{M}_{2,k}^N \Rightarrow M_{2,k}$ as $N \to \infty$, which completes the proof of (a).

We will next prove (b), making use of Lemma A.8. In the conditions of the lemma, we take, fixing $t$ and $k$ for the moment,

$$\beta_i(x, y) = \sum_{p=1}^{k} 1\{s_{p-1}^k < x \leq s_p^k\}(1\{t - s_p^k < x < t - x\} - (F(t-x) - F(t - s_p^k))).$$

Then

$$\bar{\beta}_i(x) = E[\beta_i(x, \eta_i)^2]$$
$$= \sum_{p=1}^{k} 1\{s_{p-1}^k < x \leq s_p^k\}(F(t - x) - F(t - s_p^k))$$
$$\times (1 - F(t - x) - F(t - s_p^k))$$

and (A.59) yields, by (A.13), (A.16) and (A.60),

$$(\text{A.63}) \quad N^{-1/2} R_m^N(t) = \tilde{M}_2^N(t) - \tilde{M}_{2,k}^N(t) \qquad \text{on } \{\hat{A}^N(t) \leq m\}.$$

By (A.63) and (A.59),

$$N^{-1}\langle R_m^N\rangle(t)$$
$$\leq N^{-1}\sum_{i=1}^{\hat{A}^N(t)}\sum_{p=1}^{k}1\{s_{p-1}^k < \hat{\tau}_i^N \leq s_p^k\}(F(t-s_{p-1}^k)-F(t-s_p^k))$$
$$= N^{-1}\sum_{p=1}^{k}(F(t-s_{p-1}^k)-F(t-s^k+p))(\hat{A}^N(s_p^k)-\hat{A}^N(s_{p-1}^k))$$
$$\leq \sup_{1\leq p\leq k}(N^{-1}\hat{A}^N(s_p^k)-N^{-1}\hat{A}^N(s_{p-1}^k)).$$

Then, by Lemma A.8, applying the Lenglart–Rebolledo inequality and (A.63), for $\eta > 0, \varepsilon > 0$,

$$P(|\tilde{M}_2^N(t)-\hat{M}_{2,k}^N(t)|>\eta)$$
$$\leq P(\hat{A}^N(t)>mN)+P(N^{-1/2}|R_m^N(t)|>\eta)$$
$$\leq P(N^{-1}\hat{A}^N(t)>m)+\frac{\varepsilon}{\eta^2}$$
$$+P\Big(\sup_{1\leq p\leq k}(N^{-1}\hat{A}^N(s_p^k)-N^{-1}\hat{A}^N(s_{p-1}^k))>\varepsilon\Big).$$

By Lemma A.5, continuity of the identity function $e(t)=t$ and the fact that $\max_{1\leq p\leq k}(s^k+p-s_{p-1}^k)\to 0$ as $k\to\infty$,

$$\lim_{m\to\infty}\limsup_{N\to\infty}P(N^1\hat{A}^N(t)>m)=0,$$
$$\lim_{k\to\infty}\limsup_{N\to\infty}P\Big(\sup_{1\leq p\leq k}(N^{-1}\hat{A}^N(s_p^k)-N^{-1}\hat{A}^N(s_{p-1}^k))>\varepsilon\Big)=0,$$

ending the proof of (b).

We next prove part (c). The proof proceeds in a similar manner to the proof of parts (a) and (b). Letting $B^N(t)=\lfloor Nt\rfloor$, we first note that

$$\hat{M}_2^N(t)=N^{-1/2}\int_0^t\int_0^t 1\{s+x\leq t\}\,dU^N(B^N(s),x).$$

Furthermore, setting $\bar{B}^N=\{N^{-1}B^N(t),t\geq 0\}$, it is clear that

(A.64) $$\bar{B}^N\Rightarrow e\quad\text{as }N\to\infty.$$

Next, letting

$$\check{B}_{2,k}^N(t)=\sum_{i=1}^{k}\square\tilde{U}^N((\hat{B}^N(s_{i-1}^k),0),(B^N(s_i^k),F(t-s_i^k))),$$

it follows by (A.64), Lemma 3.1 of [14] and the continuity of $\tilde{U}$ and $e$, that by (A.60),

$$(A.65) \quad \lim_{N \to \infty} P\left(\sup_{0 \leq t \leq T} |\check{B}_{2,k}^N(t) - \check{M}_{2,k}^N(t)| > \varepsilon\right) = 0, \qquad T > 0, \varepsilon > 0.$$

A similar proof to that of part (b) above can also be used to show that

$$(A.66) \quad \lim_{k \to \infty} \limsup_{N \to \infty} P(|\check{B}_{2,k}^N(t) - \hat{M}_2^N(t)| > \eta) = 0 \qquad \text{for } \eta > 0, t > 0.$$

Part (b), (A.62), (A.65) and (A.66) above now imply part (c).

Parts (a)–(c) imply the finite-dimensional convergence $(\tilde{M}_2^N, \hat{M}_2^N) \Rightarrow^{df} (\tilde{M}_2, \tilde{M}_2)$ as $N \to \infty$. It therefore remains to show that the sequence $\{(\tilde{M}_2^N, \hat{M}_2^N)\}$ is tight in order to complete the proof. However, by Proposition A.1, the sequence $\{\tilde{M}_2^N\}$ is tight and a similar if not identical proof also shows that $\{\hat{M}_2^N\}$ is tight. Thus, the sequence $\{(\tilde{M}_2^N, \hat{M}_2^N)\}$ is tight in $(D^2[0,\infty), d_{J_1}^2)$, which completes the proof. $\square$

The remainder of the appendix is now devoted to providing proofs of Lemmas A.9 and A.10. We begin with Lemma A.10. Recall first the definition of $f : \mathbb{R} \times D^3[0,\infty) \mapsto \mathbb{R} \times D^3[0,\infty)$ in (5.18) and (5.19) as

$$(A.67) \qquad f((x_1, x_2, x_3, x_4)) = (f_1(x_1), f_2(x_2), f_3(x_3), f_4(x_4)),$$

for $(x_1, x_2, x_3, x_4) \in \mathbb{R} \times D^3[0,\infty)$, where $f_1(x_1) = x_1$, $f_2(x_2) = x_2$, $f_4(x_4) = x_4$ and

$$(A.68) \qquad f_3(x_3)(\cdot) = \int_0^\cdot G(\cdot - s) \, dx_3(s),$$

where the above integral above may be interpreted as the result of integration by parts.

We then have the following result.

LEMMA A.9. *The function $f$ defined by (A.67) and (A.68) is measurable as a map from $(\mathbb{R} \times D^3[0,\infty), \mathcal{B}(\mathbb{R}) \times \mathcal{D}^3)$ to $(\mathbb{R} \times D^3[0,\infty), \mathcal{B}(\mathbb{R}) \times \mathcal{D}^3)$. Furthermore, it is continuous at continuous limits points $(x_1, x_2, x_3, x_4) \in \mathbb{R} \times D^3[0,\infty)$ such that $x_2, x_3, x_4 \in C[0,\infty)$.*

PROOF. We first show that the function $f : (\mathbb{R} \times D^3[0,\infty), \mathcal{B}(\mathbb{R}) \times \mathcal{D}^3) \mapsto (\mathbb{R} \times D^3[0,\infty), \mathcal{B}(\mathbb{R}) \times \mathcal{D}^3)$ is measurable. It is clear that $f_1 : (\mathbb{R}, \mathcal{B}(\mathbb{R})) \mapsto (\mathbb{R}, \mathcal{B}(\mathbb{R}))$, $f_2 : (D[0,\infty), \mathcal{D}) \mapsto (D[0,\infty), \mathcal{D})$ and $f_4 : (D[0,\infty), \mathcal{D}) \mapsto (D[0,\infty), \mathcal{D})$ are measurable since each of these functions are the identity functions. Therefore, if we may now show that $f_3 : (D[0,\infty), \mathcal{D}) \mapsto (D[0,\infty), \mathcal{D})$ is measurable, then, by (A.67), we will have shown the measurability of $f$.

In order to show that $f_3 : (D[0,\infty), \mathcal{D}) \mapsto (D[0,\infty), \mathcal{D})$ is measurable, first note that, integrating by parts, we have

$$f_3(x_3) = a(x_3) - b(x_3) - c(x_3), \tag{A.69}$$

where $a(x_3) = G(0)x_3, b(x_3) = x(0)G$ and

$$c(x_3)(\cdot) = \int_0^\cdot x_3(\cdot - s)\, dG(s). \tag{A.70}$$

The functions $a(x_3)$ is measurable since it is the identity function multiplied by a constant and $b(x_3)$ is measurable as well since by [1] we have that the projection map is $\pi_0(x_3) = x_3(0)$ is a measurable function too. Thus, since the sum of a finit number of measurable functions is measurable, it remains to show that $c(x_3)$ is measurable from $(D[0,\infty), \mathcal{D})$ to $(D[0,\infty), \mathcal{D})$ in order to complete the proof of the measurability of $f$. However, this may be shown in a manner similar to the proof of the measurability of $\Psi_B^a$ in the measurability portion of the proof of Proposition 3.1. We omit the details for the sake of brevity. This completes the proof of the measurability of $f$.

We now show that $f : (\mathbb{R} \times D^3[0,\infty), |\cdot| \times d_{J_1}^3) \mapsto (\mathbb{R} \times D^3[0,\infty), |\cdot| \times d_{J_1}^3)$ is continuous at continuous limit points $(x_1, x_2, x_3, x_4) \in \mathbb{R} \times D^3[0,\infty)$ such that $x_2, x_3, x_4 \in C[0,\infty)$. First, it is clear that the functions $f_1 : (\mathbb{R}, |\cdot|) \mapsto (\mathbb{R}, |\cdot|), f_2 : (D[0,\infty), u) \mapsto (D[0,\infty), u)$ and $f_4 : (D[0,\infty), u) \mapsto (D[0,\infty), u)$ are continuous. This follows easily since each of these functions are the identity functions. We now show that the function $f_3 : (D[0,\infty), u) \mapsto (D[0,\infty), u)$ is continuous. This then implies that $f : (\mathbb{R} \times D^3[0,\infty), |\cdot| \times u^3) \mapsto (\mathbb{R} \times D^3[0,\infty), |\cdot| \times u^3)$ is continuous. However, since converge in $(D[0,\infty), d_{J_1})$ to a continuous limit point $x \in C[0,\infty)$ is equivalent to convergence in $(D[0,\infty), u)$, this then implies that $f : (\mathbb{R} \times D^3[0,\infty), |\cdot| \times u^3) \mapsto (\mathbb{R} \times D^3[0,\infty), |\cdot| \times d_{J_1}^3)$ is continuous a continuous limit points $(x_1, x_2, x_3, x_4) \in \mathbb{R} \times C^3[0,\infty)$, which completes the proof.

Suppose first that $x_3^n \to x_3$ in $(D[0,\infty), u)$ as $n \to \infty$. It then follows that for each $T \geq 0$,

$$\sup_{0 \leq t \leq T} |x_3^n(t) - x(t)| \to 0 \quad \text{as } n \to \infty.$$

However, recalling the representation of $f_3(x_3)$ in (A.70), this then implies that

$$\sup_{0 \leq t \leq T} |f_3(x_3^n)(t) - f_3(x_3)(t)|$$
$$= \sup_{0 \leq t \leq T} |(a(x_3^n)(t) - b(x_3^n)(t) - c(x_3^n)(t))$$
$$\qquad - (a(x_3)(t) - b(x_3)(t) - c(x_3)(t))|$$
$$= \sup_{0 \leq t \leq T} |(a(x_3^n)(t) - a(x_3)(t)) - (b(x_3^n)(t) - b(x_3)(t))$$

$$- (c(x_3^n)(t) - c(x_3)(t))|$$

$$= \sup_{0 \leq t \leq T} \Big| (G(0)x_3^n(t) - G(0)x_3(t)) - (x_3^n(0)G(t) - x_3(0)G(t))$$

$$- \Big( \int_0^t x_3^n(t-s)\,dG(s) - \int_0^t x_3(t-s)\,dG(s) \Big) \Big|$$

$$= \sup_{0 \leq t \leq T} \Big| G(0)(x_3^n(t) - x_3(t)) - G(t)(x_3^n(0) - x_3(0))$$

$$- \int_0^t (x_3^n(t-s) - x_3(t-s))\,dG(s) \Big|$$

$$\leq \sup_{0 \leq t \leq T} |G(0)(x_3^n(t) - x_3(t))| + \sup_{0 \leq t \leq T} |G(t)(x_3^n(0) - x_3(0))|$$

$$+ \sup_{0 \leq t \leq T} \Big| \int_0^t (x_3^n(t-s) - x_3(t-s))\,dG(s) \Big|$$

$$\leq \sup_{0 \leq t \leq T} |G(0)(x_3^n(t) - x_3(t))| + \sup_{0 \leq t \leq T} |G(t)(x_3^n(0) - x_3(0))|$$

$$+ \sup_{0 \leq t \leq T} \int_0^t |x_3^n(t-s) - x_3(t-s)|\,dG(s)$$

$$\leq G(0) \sup_{0 \leq t \leq T} |x_3^n(t) - x_3(t)| + G(T)|x_3^n(0) - x_3(0)|$$

$$+ G(T) \sup_{0 \leq t \leq T} |x_3^n(t) - x_3(t)|$$

$$= (G(0) + 2G(T)) \sup_{0 \leq t \leq T} |x_3^n(t) - x_3(t)|$$

$$\to 0 \quad \text{as } n \to \infty,$$

and so the function $f_3$ is continuous as a map from $(D[0,\infty), u)$ to $(D[0,\infty), u)$, which completes the proof. $\square$

We now provide a proof of Lemma A.10 which is instrumental in the proof of Proposition 5.3. Our setup is the same as in the proof of Lemma A2 of Puhalskii and Reiman [17]. In particular, we consider a sequence of queueing systems indexed by $N$, each operating under the FIFO service discipline and each with a single arrival process $A^N = \{A^N(t), t \geq 0\}$ and a single departure process $D^N = \{D^N(t), t \geq 0\}$. We denote by $Q^N(t)$ the queue length of the $N$th system at time $t$, by $V^N(t)$ the virtual waiting time at time $t$, and by $V_i^N$ the waiting time of the $i$th customer to arrive to the system. Finally, we set $\tilde{Q}^N(t) = N^{-1/2}Q^N(t)$, $\tilde{A}^N(t) = N^{-1/2}(A^N(t) - \lambda^N t)$ where $\{\lambda^N\}$ is a sequence of constants and we assume that $A^N(0) = D^N(0) = 0$. We then have the following result.

LEMMA A.10. *If $\tilde{Q}^N \Rightarrow \tilde{Q}$ and $\tilde{A}^N \Rightarrow \tilde{\xi}$ as $N \to \infty$, each on $D[0, \infty)$, and $\lambda^N/N \to \lambda > 0$, then the processes $\{N^{1/2}V^N(t), t \geq 0\}$ and $\{N^{1/2}V^N_{\lfloor Nt \rfloor}, t \geq 0\}$ converge in distribution on $(D[0, \infty), d_{J_1})$ to the respective processes $\{\tilde{Q}(t)/\lambda, t \geq 0\}$ and $\{\tilde{Q}(t/\lambda)/\lambda, t \geq 0\}$.*

PROOF. First, note that since $\tilde{Q}^N \Rightarrow \tilde{Q}$ and $\tilde{A}^N \Rightarrow \tilde{\xi}$ as $N \to \infty$, it follows by Prohorov's theorem [1] that the sequences $\{\tilde{Q}^N, N \geq 1\}$ and $\{\tilde{A}^N, N \geq 1\}$ are both tight. Thus, we have that sequence $\{(\tilde{Q}^N, \tilde{A}^N), N \geq 1\}$ is tight in $(D^2[0, \infty), d^2_{J_1})$ and hence, by a second application of Prohorov's theorem [1], the sequence $\{(\tilde{Q}^N, \tilde{A}^N), N \geq 1\}$ is relatively compact. Thus, for any subsequence $\{N_k\}$, there exists a further subsequence $\{N'_k\}$ such that $(\tilde{Q}^{N'_k}, \tilde{A}^{N'_k}) \Rightarrow (\hat{\tilde{Q}}, \hat{\tilde{A}})$ as $k \to \infty$. Thus, by Lemma A2 of Puhalskii and Reiman [17], we have that $\{N'^{1/2}_k V^{N'_k}(t), t \geq 0\}$ and $\{N'^{1/2}_k V^{N'_k}_{\lfloor Nt \rfloor}, t \geq 0\}$ converge in distribution on $(D[0, \infty), d_{J_1})$ to the respective processes $\{\hat{\tilde{Q}}(t)/\lambda, t \geq 0\}$ and $\{\hat{\tilde{Q}}(t/\lambda)/\lambda, t \geq 0\}$. However, since it must be the case that $\hat{\tilde{Q}} \stackrel{d}{=} \tilde{Q}$ and the sequence $\{N_k\}$ was arbitrary, this completes the proof. □

**Acknowledgments.** This work was part of the author's Ph.D. thesis. The author would therefore like to thank his Ph.D. advisors, Jim Dai and Amy Ward, for their constant guidance, support and encouragement. This work would surely not have been possible without either one of them. In addition, the author would also like to thank Bob Foley, Marty Reiman and Ward Whitt for their numerous helpful comments and suggestions on earlier versions of this draft and the two anonymous referees whose guidance greatly improved the exposition of the paper.

LEONARD N. STERN SCHOOL OF BUSINESS
NEW YORK UNIVERSITY
HENRY KAUFMAN MANAGEMENT CENTER
44 WEST FOURTH STREET
NEW YORK, NEW YORK 10012
E-MAIL: jreed@stern.nyu.edu